\documentclass{article}

\usepackage{delarray,verbatim,enumerate,a4wide}
\usepackage{amsmath,amsthm,amstext,amsbsy,amssymb,amsfonts,amscd}
\usepackage[all]{xy}
\usepackage[frenchb]{babel}
\usepackage[T1]{fontenc}
\usepackage[utf8]{inputenc}
\usepackage{color}
\definecolor{vert}{rgb}{0.1,0.4,0.2}
\usepackage[colorlinks=true,linkcolor=blue,citecolor=vert]{hyperref}

\usepackage{calligra}
\DeclareFontShape{T1}{calligra}{m}{n}{<->s*[0.95]callig15}{}
\DeclareMathAlphabet{\mathscr}{T1}{calligra}{m}{n}

\newtheorem{Th}{Théorème}[]

\newtheorem{Prop}[Th]{Proposition}
\newtheorem{Cor}[Th]{Corollaire}
\newtheorem{Conj}[Th]{Conjecture}
\newtheorem{Sco}[Th]{Scolie}
\newtheorem{Def} [Th]{Définition}
\newtheorem{PDef}[Th]{Proposition \& Définition}
\newtheorem{DProp}[Th]{Définition \& Proposition}

\newtheorem{DTh}[Th]{Définition \& Théorème}
\newtheorem{ThD}[Th]{Théorème \& Définition}
\newtheorem{Exe}[Th]{Exemple}
\newtheorem*{Th*}{Théorème} 
\newtheorem*{Sco*}{Scolie}

\def\Preuve{\smallskip\noindent {\it Preuve.~}}

\def\Remarque{\smallskip\noindent {\it Remarque.~}}

\font\teneufm=eufm10
\font\seveneufm=eufm7
\font\fiveeufm=eufm5
\newfam\gothfam
\textfont\gothfam=\teneufm \scriptfont\gothfam=\seveneufm
\scriptscriptfont\gothfam=\fiveeufm

		\def\QQ{\mathbb Q}	
\def\NN{\mathbb N}	\def\ZZ{\mathbb Z}
\def\F2{\mathbb{F}_2}	\def\Z2{\mathbb{Z}_2}	\def\lT{\ov{\mathbb T}_\ell}	
\def\Zl{\mathbb{Z}_\ell} 	\def\Ql{\mathbb{Q}_\ell}	\def\Tl{\mathbb{T}_\ell}

 				\def\U{\mathcal  U}	\def\F{\mathcal  F}
\def\J{\mathcal  J}  	\def\C{\mathcal  C}	\def\R{\mathcal  R}	\def\N{\mathcal N}
 	\def\Pl{\mathcal  P\ell}  	\def\Cl{\mathcal  C \ell}
\def\E{\mathcal  E}		\def\T{\mathcal  T}		

	\def\p{{\mathfrak p}}	\def\x{{\mathfrak x}}		\def\y{{\mathfrak y}}
		\def\l{{\mathfrak l}}	

\def\W{W\!K}

\def\wi{\widetilde}	\def\ov{\overline}	
\def\rg{\operatorname{rg}}	\def\div{\operatorname{div}}
	\def\deg{\operatorname{deg}}
\def\Gal{\operatorname{Gal}}	\def\Log{\operatorname{Log}}	\def\Rad{\operatorname{Rad}}
\def\Ker{\operatorname{Ker}}	\def\Coker{\operatorname{Coker}}	\def\Hom{\operatorname{Hom}}

\makeatletter
\newcommand*\wt[2][0.2ex]{%
        \begingroup
        \mathchoice{\wt@helper{#1}{#2}{\displaystyle}{\textfont}}
                   {\wt@helper{#1}{#2}{\textstyle}{\textfont}}
                   {\wt@helper{#1}{#2}{\scriptstyle}{\scriptfont}}
                   {\wt@helper{#1}{#2}{\scriptscriptstyle}{\scriptscriptfont}}%
        \endgroup
        #2%
}
\newcommand*\wt@helper[4]{%
        \def\currentfont{\the#41}%
        \def\currentskewchar{\char\the\skewchar\currentfont}%
        \setbox\tw@\hbox{\currentfont$#2$\currentskewchar}%
        \dimen@ii\wd\tw@
        \setbox\tw@\hbox{\currentfont$#2${}\currentskewchar}%
        \advance\dimen@ii-\wd\tw@
        \rlap{\raisebox{-#1}{$\m@th#3\kern\dimen@ii\widetilde{\phantom{#2}}$}}%
}
\makeatother

\def\wE{\,\wt[0.1ex]{\!\mathcal E}}	\def\we{\wt[0.3ex]{\mathfrak E}}	\def\wU{\wt[0.2ex]{\mathcal U}}
\def\wJ{\,\wt[0.2ex]{\!\mathcal J}}	\def\wCl{\wt[0.1ex]{\mathcal C\!\ell}} \def\wDl{\wt[0.2ex]{\mathcal D\!\ell}}
\def\bot{{\scriptscriptstyle\perp}}

\begin{document}

\title{\LARGE\bf Sur les normes cyclotomiques et les conjectures de Leopoldt et de Gross-Kuz'min}

\author{ Jean-François {\sc Jaulent} }
\date{}
\maketitle
\bigskip\bigskip\bigskip

{\small
\noindent{\bf Résumé.} À la lumière de la Théorie $\ell$-adique du corps de classes, nous dressons un bref panorama des propriétés arithmétiques des groupes de normes cyclotomiques en liaison avec les conjectures de Leopoldt et de Gross. Nous revenons pour les compléter sur quelques résultats classiques que nous illustrons au moyen des classes et unités logarithmiques par divers exemples et contre-exemples obtenus avec PARI.
}

\

{\small
\noindent{\bf Abstract.} We use $\ell$-adic class field theory to take a new view on cyclotomic norms and Leopoldt or Gross-Kuz'min conjectures. By the way we recall and complete some classical results. We illustrate the logarithmic approach by various numerical examples and counter-examples obtained with PARI.
}
\bigskip\bigskip

\noindent{\large \bf Introduction}
\medskip

Le groupe des $\ell$-normes cyclotomiques, i.e. des éléments du $\ell$-adifié $\R_K=\Zl\otimes_\ZZ K^\times$ du groupe multiplicatif d'un corps de nombres $K$ qui sont normes dans tous les étages finis $K_n/K$ de la $\Zl$-extension cyclotomique $K_\infty$ de ce corps, joue un rôle crucial dans notre compréhension de l'arithmétique de ses pro-$\ell$-extensions abéliennes.\smallskip

Or ce groupe est de nature irréductiblement $\ell$-adique et ne peut être compris en profondeur si on tente de le  le regarder (ou si l'on essaie de le définir) comme le $\ell$-adifié d'un sous-groupe  naturel (ou naïf) du groupe multiplicatif $K^\times$: c'est intrinsèquement un sous-module de $\R_K$ qui ne provient pas en général (par $\ell$-adification) d'un tel sous-groupe du groupe $K^\times$.\smallskip

La Théorie $\ell$-adique du corps de classes (cf. \cite{J18,J31} ou \cite{Gra}), qui permet de travailler avec les objets $\ell$-adiques comme s'ils étaient des objets classiques, fournit en revanche un cadre conceptuel efficace pour en mener l'étude. Plus précisément, l'introduction des classes et unités logarithmiques permet d'interpréter ce groupe des normes cyclotomiques comme un groupe d'unités, analogue au $\ell$-adifié $\E_K=\Zl\otimes_\ZZ E_K$ du groupe des unités au sens habituel, en lien avec un groupe de classes, qui est le pendant logarithmique du $\ell$-groupe des classes d'idéaux au sens habituel.\smallskip

Ces groupes logarithmiques possèdent un comportement à certains égards très voisin de ceux des groupes classiques (cf. \cite{J28}), quelquefois différents (cf. \cite{J53}). Ils mettent en jeu des conjectures encore ouvertes -- on ne sait pas, par exemple si le $\ell$-groupe des classes logarithmiques est fini en dehors du cas abélien et de quelques autres (cf. \cite{J10,J17}) -- mais sont totalement effectifs en ce sens que pour un corps de nombres $K$ et un premier $\ell$ donnés, on dispose d'algorithmes permettant de les calculer tout comme leurs analogues standard. De plus, comme observé dans \cite{Grt}, de par leurs propriétés de descente, ces groupes se révèlent bien adaptés à la théorie d'Iwasawa cyclotomique.\smallskip

Le but de cette note est ainsi de revisiter un certain nombre de résultats connus, méconnus ou carrément ignorés, à la lumière de l'arithmétique logarithmique qui apporte des informations précises sur des objets a priori définis en haut de la $\Zl$-extension cyclotomique $K_\infty$, mais que l'on peut ainsi (au moins partiellement) appréhender à l'aide d'invariants arithmétiques ne faisant intervenir que le corps $K$, et d'illustrer le tout par des exemples numériques éclairants qui contredisent en particulier une conjecture de Seo \cite{Se1}.
L'approche retenue peut ainsi être regardée comme complémentaire de celle de Movahhedi et Nguyen Quang Do \cite {MN2} sur les mêmes sujets.


\newpage
\noindent{\large \bf 1. Normes cyclotomiques et unités logarithmiques}
\medskip

Soient $K$ un corps de nombres et $K_\infty=\cup_{n\in\NN}K_n$ sa $\Zl$-extension cyclotomique. Notons $\Gamma=\gamma^{\Zl}$ le groupe de Galois $\Gal(K_\infty/K)$, identifié à $\Zl$ par le choix d'un générateur topologique $\gamma$ et $\Gamma_{\!n}=\Gamma^{\ell^n}$ le sous-groupe ouvert $\Gal(K_\infty/K_n)$.

 Pour appréhender les problèmes normiques dans la pro-$\ell$-extension $K_\infty/K$, nous nous appuierons sur la Théorie $\ell$-adique du corps de classes, telle qu'exposée dans \cite{J18} et \cite{J31}, qui permet de raisonner sur les extensions infinies comme sur celles finies.\smallskip

Notons donc $\J _K$ le {\em $\ell$-adifié  du groupe  des idéles de $K$}, i.e. le produit $\J_K=\prod_\p ^{res}\R_{K_\p}$ des compactifés $\ell$-adiques respectifs $\R _{K_\p}=\varprojlim K_\p ^\times /K_\p^{\times \ell^m}$ des groupes multiplicatifs des complétés $K_\p$, restreint aux familles $(x_\p)_\p$ dont presque tous les éléments tombent dans le sous-groupe unité.

\begin{DProp}
{\rm Pour chaque place finie $\p$ le sous-groupe $\wU _{K_\p}$ de $\R _{K_\p }$ formé des normes cyclotomiques locales, i.e. des éléments de $\R _{K_\p}$  qui sont normes à chaque étage fini de la $\Zl$-extension cyclotomique locale $K_\p ^c/K_\p$,
est par convention le} groupe des unités logarithmiques  {\rm de $K_\p$ car c'est le noyau de la valuation logarithmique $\wi\nu_\p$ à valeurs dans $\Zl$ définie par la formule:\smallskip

\centerline{$\tilde\nu_\p :  \, x_\p\, \mapsto\, -\frac{\Log_\ell\, |  x_\p | _\p}{\deg\, \p} =  -\frac{\Log_\ell\, N_{K_\p/\QQ_p}(x_\p)}{\deg\, \p}$,}\smallskip

\noindent qui fait intervenir le logarithme d'Iwasawa de la norme locale de $x_\p$ (cf. \cite{J28,J31}).}
\end{DProp}

La Théorie $\ell$-adique du corps de classes établit un isomorphisme de groupes topologiques compacts entre le $\ell$-groupe des classes d'idèles $\,\C_K$ défini comme quotient\smallskip

\centerline{$\C_K=\J_K/\R_K$}\smallskip

\noindent de $\J_K$ par son sous-groupe principal $\R_K=\Zl\otimes_\ZZ K^\times$ et le groupe de Galois $G_K=\Gal(K^{ab}/K)$ de la pro-$\ell$-extension abélienne maximale de $K$. Ainsi  (cf. \cite{J28,J31}):

\begin{ThD}\label{Caractérisation}
Dans la correspondance du corps de classes $\ell$-adique:
\begin{itemize}
\item[(i)] le groupe de normes associé à la plus grande sous-extension $K^{lc}$ de $K^{ab}$ qui est localement cyclotomique (i.e. complètement décomposée sur $K_\infty$ en chacune de ses places) est le produit $\,\wU_K\R_K$ du sous-groupe $\,\wU_K=\prod_\p \, \wU _{K_\p}$ des unités logarithmiques locales et de $\R_K$;
\item[(ii)] le groupe de normes associé à la $\Zl$-extension cyclotomique $K^c=K_\infty$ est le sous-groupe des idèles de degré nul: $\wJ_K=\{\x=(x_\p)_\p\in\J_K\,|\, \deg(\x)=\sum_\p\wi\nu(x_\p)\deg\,\p=0\}$.
\end{itemize}
En particulier, le groupe de Galois $\Gal(K^{lc}/K^c)$ s'identifie au quotient $\,\wCl_K= \wJ_K/\wU_K\R_K$. 

Nous disons que $\,\wCl_K$ est le $\ell$-groupe des classes logarithmiques du corps $K$.
\end{ThD}

 Le groupe $\wCl_K$ peut être regardé comme le quotient du groupe $\wDl _K = \wJ _K /\wU_K$ des diviseurs logarithmiques (de degré nul) par son sous-groupe principal $\Pl_K=\R_K\wU_K/\wU_K$; le numérateur $\,\wDl _K$  s'identifie au sous-groupe $\wi\oplus_\p \,\Zl\,\p$ des diviseurs de degré nul de la somme formelle $\oplus_\p \,\Zl\,\p$.

\begin{DTh}
Le noyau $\wE _K=\R_K \cap\, \wU_K$ du morphisme $\wi\div:\;x\mapsto \sum_\p\wi\nu_\p(x)\,\p$ de $\R_K$ dans $\wDl_F$ est le groupe des unités  logarithmiques  (globales) du corps $K$. Il est donné par la suite exacte:

\centerline{$1 \rightarrow \wE_K \rightarrow \R_K \rightarrow \wDl_K=\wi\oplus_\p \,\Zl\,\p \rightarrow \wCl_K \rightarrow 1$.}\smallskip

Le groupe $\wE_K$ est exactement le sous-groupe des normes cyclotomiques globales de $\R_K$.
\end{DTh}

\noindent{\em Preuve.} Par construction les éléments de $\wE_K$ sont, en effet, les éléments globaux qui sont normes locales partout dans chacune des sous-extensions finies $K_n/K$ de l'extension procyclique $K_\infty/K$, donc normes globales dans $K_\infty/K$, en vertu du principe $\ell$-adique de Hasse (cf. Th. \ref{PH} ci-après).\smallskip

\Remarque Plusieurs auteurs (ainsi \cite{Grt,Kat,Kuz,MN2,MR}) utilisent l'appellation {\em normes universelles} pour désigner les éléments qui sont normes à chaque étage fini d'une $\Zl$-extension donnée, et même plus spécifiquement normes (locales ou globales) de $\ell$-unités: ce sont, en effet, les normes universelles pour un certain foncteur  (cf. \S 7). Il convient de ne pas les confondre avec les {\em normes universelles} au sens de la Théorie du corps de classes: dans la théorie classique (cf. e.g. \cite{AT}), ces normes universelles sont les éléments de la composante connexe du neutre dans le groupe des classes d'idèles; dans la théorie $\ell$-adique  (cf. \cite{Gra,J18,J31}), cette composante se réduit au neutre.

\newpage
\noindent{\large \bf2.  La conjecture de Gross-Kuz'min}
\medskip

\begin{Def}[{\bf Module de Kuz'min-Tate}]
Nous appelons module de Kuz'min-Tate d'un corps de nombres $K$  le groupe de Galois $\,\T_K=\Gal(K_\infty^{cd}/K_\infty)$ de la pro-$\ell$-extension abélienne maximale de la $\Zl$-extension cyclotomique $K_\infty$ de $K$ qui est complètement décomposée en toutes les places, groupe que  la Théorie du corps de classes identifie à la limite projective pour les applications normes $\,\varprojlim \Cl'_{K_n}$ des $\ell$-groupes de $\ell$-classes $\Cl'_{K_n}$ attachés aux étages finis $K_n$ de la tour $K_\infty/K$. 
\end{Def}

L'appellation {\em module de Tate} pour désigner le groupe $\T_K$ a été introduite par Kuz'min \cite{Kuz}. Mais d'autres traditions existent dans le cadre de la Théorie d'Iwasawa. Ainsi dans \cite{J18,J23}, elle est utilisée pour désigner la limite projective $\mathbb T_\ell=\varprojlim\mu_{K_n}^{\phantom{loc}}$ des $\ell$-groupes de racines de l'unité dans la tour cyclotomique $K_\infty/K$, lorsque le corps $K$ contient les racines 2$\ell$-ièmes de l'unité. C'est pourquoi, pour prévenir toute confusion, nous parlons ici de {\em module de Kuz'min-Tate}.\smallskip

Ce point de vocabulaire précisé, nous avons (cf. \cite{FG,Gro,J18,J28,J31,Ko1,Kuz}):

\begin{Th}[{\bf Conjecture de Gross-Kuz'min}]\label{CGK}
Pour tout corps de nombres $K$ et tout premier $\ell$, les conditions suivantes sont équivalentes:
\begin{itemize}
\item[(i)] Le polynôme caractéristique $\chi_{\T_K}(X)$ du $\Lambda$-module $\T_K$ n'est pas divisible par $X$: $\chi_{\T_K}(0)\ne 0$.
\item[(ii)] Le sous-groupe ambige $\T_K^{\,\Gamma}$ est fini.$\phantom{\wCl}$
\item[(iii)] Le quotient des genres ${}_\Gamma\T_K = \T_K/\T_K^{\,(\gamma-1)}\simeq \, \wCl_K$ est fini.
\end{itemize}
Lorsqu'elles sont réalisées, i.e. lorsque le $\ell$-groupe des classes logarithmiques $\;\wCl_K$ est fini, nous disons que le corps $K$ satisfait la conjecture de Gross-Kuz'min  pour le premier $\ell$.
\end{Th}

\noindent{\em Preuve.} L'équivalence résulte de la théorie d'Iwasawa; l'isomorphisme  ${}_\Gamma\T_K\simeq\wCl_K$ de l'interprétation galoisienne du groupe  des classes logarithmiques donnée par le Théorème \ref{Caractérisation}.\medskip

\noindent {\em Nota.} Cette conjecture est référée sous diverses appellations suivant le contexte comme conjecture de Kuz'min, de Gross ou encore de Gross généralisée. On sait notamment par des arguments de transcendance qu'elle est vérifiée par les corps absolument abéliens et quelques autres.

\begin{Sco}\label{Ulog}
Le groupe $\,\wE_K$ des unités logarithmiques de $K$ est le produit:

\centerline{ $\wE_K= \mu^{\phantom{loc}}_K  \Zl^{r_K+c_K+\delta_K^{\mathscr G}}$}\smallskip

\noindent du $\ell$-sous-groupe $\mu^{\phantom{lc}}_K$ des racines globales de l'unité et d'un $\Zl$-module libre de rang $r_K+c_K+\delta_K^{\mathscr G}$, où  $r_K$ et $c_K$ sont respectivement les nombres de places réelles et complexes de $K$; et $\delta_K^{\mathscr G}=\dim^{\phantom{c}}_{\Zl} \wCl_K$  mesure le défaut dans $K$ de la conjecture de Gross-Kuz'min  pour  le premier $\ell$.
\end{Sco}

\noindent{\em Preuve.} C'est une conséquence immédiate de la suite exacte restreinte aux $\ell$-unités (cf. \cite{J28,J31}):\smallskip

\centerline{$1 \rightarrow \wE_K \rightarrow \E'_K \rightarrow \wDl_K=\wi\oplus_{\p|\ell} \,\Zl\,\p \rightarrow \wCl_K \rightarrow \Cl'_K$.}
\medskip

\Remarque La plupart des modules $\ell$-adiques globaux que nous manipulons sont les $\ell$-adifiés des groupes usuels attachés à un corps de nombres. Mais {\em ce n'est pas le cas des groupes d'unités logarithmiques}, comme le montre le contre-exemple suivant (qui contredit une affirmation de \cite{Se3}).

\begin{Prop}\label{Contrex}
L'image $\Zl\otimes_\ZZ\wi E_K$ dans $\R_K$ du sous-groupe $\wi E_K=\cap_{n\in\NN}N_{K_n/K}(K^\times_n)$ de $K^\times$ formé des éléments qui sont normes (locales comme globales) dans chaque sous-extension finie de la tour $K_\infty/K$ est contenue dans le groupe des unités logarithmiques; mais ce dernier n'en est pas, en général, le $\ell$-adifié, l'inclusion

\centerline{$\Zl\otimes_\ZZ\wi E_K \subset \wE_K$}\smallskip

\noindent pouvant être stricte; ainsi pour $K=\QQ[\sqrt d]$ corps quadratique réel dans lequel $\ell$ est décomposé.
\end{Prop}

\noindent{\em Preuve.} Soient  $K=\QQ[\sqrt d]$ un corps quadratique réel et $\ell\ge 3$  décomposé dans $K/\QQ$ . Le groupe $E'_K$ des $\ell$-unités de $K$ est le $\ZZ$-module multiplicatif engendré par les quatre éléments $-1, \varepsilon,\ell$ et $\eta$, où $\varepsilon$ est l'unité fondamentale et $\eta$ un générateur de la plus petite puissance principale non triviale de l'une des places au-dessus de $\ell$. Notons $\bar\eta$ le conjugué de $\eta$. Pour $(a,b,c)\in\ZZ^3$, l'implication\smallskip

\centerline{$\Log_\ell(\pm \varepsilon^a \ell^b\eta^c)=\log_\ell(\varepsilon^a/\bar\eta^c)=0 \quad \Rightarrow \quad \varepsilon^a/\bar\eta^c=1 \quad \Rightarrow \quad a=c=0$}\smallskip

\noindent montre que le sous-groupe des normes cyclotomiques {\em naïves} $\wi E_K$ est engendré par -1 et $\ell$. Il suit:\smallskip

\centerline{$\dim^{\phantom{c}}_{\ZZ}\wi E_K=1<2=\dim^{\phantom{c}}_{\Zl}\wE_K$ (puisqu'on a ici $\delta^{\mathscr G}_K=0$ par abélianité).}

\newpage
\noindent{\large \bf 3. Non trivialité du $\ell$-groupe des classes logarithmiques}
\medskip

Étudions maintenant plus attentivement la trivialité de $\wCl_K$ pour un $\ell$ fixé. Notons comme précédemment $K_\infty = \cup_{n\in\NN}K_n$ la $\Zl$-extension cyclotomique de $K$ et considérons le module de Kuz'min-Tate $\T_K=\varprojlim \Cl'_{K_n}$, c'est à dire la limite projective des $\ell$-groupes de $\ell$-classes d'idéaux des étages finis $K_n$, qui s'identifie au groupe de Galois le la pro-$\ell$-extension abélienne maximale de $K_\infty$ qui est complètement décomposée en toutes les places. 

\begin{PDef}\label{trivial} Les assertions suivantes sont équivalentes:
\begin{itemize}
\item[(i)]  Le  quotient des genres de  $\T_K$ est trivial: $\wCl_K=1$.
\item[(ii)] Les $\ell$-groupes de $\ell$-classes $\;\Cl'_{K_n}$ sont ultimement triviaux: $\Cl'_{K_n}=1$ pour $n\gg 0$.
\item[(iii)] Le module de Kuz'min-Tate $\T_K$ est trivial: $\T_K=1$.
\end{itemize}
Lorsqu'elles sont vérifiées, nous disons que le corps $K$ est $\ell$-logarithmiquement principal.$\phantom{\wCl_K}$
\end{PDef}

\noindent {\em Preuve.} La $\Zl$-extension cyclotomique étant ultimement ramifiée, l'application naturelle $\T_K\rightarrow\Cl'_{K_n}$ est ultimement surjective; d'où l'équivalence {\em (ii)} $\Leftrightarrow$ {\em (iii)}. L'équivalence {\em (i)} $\Leftrightarrow$ {\em (iii)} résulte directement, elle, du lemme de Nakayama appliqué au $\Lambda$-module noethérien $\T_K$.
\medskip

Il est très facile de vérifier qu'il existe des corps de nombres non logarithmiquement principaux: un algorithme général de calcul du groupe des $\ell$-classes logarithmiques, pour un premier $\ell$ et un corps de nombres $K$ donnés, a ainsi été développé d'abord dans \cite{DS}, puis, dans une version améliorée, dans \cite{DJ+}. Il en est de même pour $\ell=2$ de versions raffinées prenant en compte la notion de classe logarithmique signée, qui est l'analogue logarithmique des classes d'idéaux au sens restreint \cite{JPPS1, JPPS2}. Ces algorithmes sont effectifs et ne présupposent nullement la validité de la conjecture de Gross-Kuz'min : ils vérifient concrètement par le calcul même qu'elle est satisfaite pour le premier $\ell$ et le corps de nombres $K$ étudiés. Des versions ad hoc plus élémentaires peuvent être mises en \oe uvre dans les cas de tout petit degré. Ainsi:

\begin{Exe}[{\bf Corps quadratiques imaginaires non logarithmiquement principaux}]${}$\par

\noindent Nous listons ci-dessous quelques exemples de couples $(K,\ell)$ pour lesquels $\wCl_K$ est non trivial.
\begin{itemize}
\item[(i) \;]    Pour $K=\QQ [\sqrt{-1}]$, le plus petit premier pour lequel $\;\wCl_K$ est non trivial est: $\ell= 29\, 789$.
\item[(ii) ]   Pour $K=\QQ [\sqrt{-5}]$, le plus petit premier pour lequel $\;\wCl_K$ est non trivial est: $\ell= 5\; 881$.
\item[(iii)]   Pour $K=\QQ [\sqrt{-7}]$, le plus petit premier pour lequel $\;\wCl_K$ est non trivial est: $\ell= 19\; 531$.
\item[(iv)]  Pour $K=\QQ [\sqrt{-11}]$, le plus petit premier pour lequel $\;\wCl_K$ est non trivial est: $\ell= 5$.
\item[(v) ]  Pour $K=\QQ [\sqrt{-13}]$, le plus petit premier pour lequel $\;\wCl_K$ est non trivial est: $\ell= 113$.
\item[(vi)] Pour $K=\QQ [\sqrt{-31}]$, le plus petit premier pour lequel $\;\wCl_K$ est non trivial est: $\ell= 227$.
\end{itemize}
Pour chacun d'eux, le $\ell$-groupe des classes logarithmiques est d'ordre exactement $\ell$.
\end{Exe}

Les calculs --très simples-- conduisant à ces résultats ont été effectués à Bordeaux par Bill Allombert à l'aide du logiciel PARI. Comme on peut le voir, la taille du plus petit $\ell$ donnant un groupe non trivial n'est pas corrélée simplement au discriminant du corps considéré. Les heuristiques amènent par ailleurs à penser que, pour chaque corps quadratique imaginaire $K$, il existe en fait une infinité de premiers $\ell$ pour lesquels le $\ell$-groupe $\wCl_K$ est non trivial. Ainsi:

\begin{Exe}[{\bf Groupes de classes logarithmiques du corps quadratique $\QQ[\sqrt{-3}]$}]   Pour $K=\QQ [\sqrt{-3}]$, les plus petits premiers $\ell$ pour lequel $\;\wCl_K$ est non trivial sont: $13$, $181$, $2\,521$, $76\,543$, $489\,061$. Dans tous ces cas le $\ell$-groupe des classes logarithmiques est d'ordre exactement $\ell$.
\end{Exe}

Pour les corps réels, en revanche, on peut toutefois risquer une conjecture différente:

\begin{Conj}
Pour chaque corps quadratique réel $K$ donné, il existe un nombre fini de premiers $\ell$ pour lesquels le $\ell$-groupe des classes logarithmiques de $K$ est non trivial.
\end{Conj}

Comme un tel corps vérifie trivialement la conjecture de Gross-Kuz'min, la conjecture ci-dessus postule donc que le groupe des classes logarithmiques {\em global} d'un corps quadratique réel, i.e. le produit des $\ell$-groupes de classes logarithmiques  pour tous les premiers $\ell$, est un groupe fini.\medskip

\newpage
\noindent{\large \bf 4. Énoncés de la conjecture de Leopoldt}
\medskip

La conjecture initiale de Leopoldt a été énoncée pour les corps totalement réels. Mais elle vaut, en fait, mutatis mutandis, quelle que soit la signature du corps de nombres considéré. Elle a alors une traduction particulièrement élégante en termes de corps de classes.

Rappelons que la Théorie $\ell$-adique du corps de classes (cf. \cite{J18,J31}) assure que le $\ell$-adifié $\R_K=\Zl\otimes K^\times$ du groupe multiplicatif d'un corps de nombres $K$ s'injecte naturellement dans le $\ell$-adifié $\J_K$ du groupe des idèles (défini comme le produit restreint $\J_K=\prod^{res}_\p\R_{K_\p}$ des compactifiés $\ell$-adiques $\R_{K_\p}=\varprojlim K^\times_\p/K^{\times \ell^m}_\p$  des groupes multiplicatifs respectifs des complétés de $K$), ce qui permet de considérer l'intersection\smallskip

\centerline{$\mu_K^{loc}=\R_K\cap\prod_\p\mu^{\phantom{loc}}_{K_\p}$}\smallskip

\noindent de $\R_K$ avec les produit des sous-groupes locaux de racines de l'unités (dans les compactifiés $\R_{K_\p}$), c'est-à-dire le sous-groupe des éléments {\em globaux} qui sont {\em localement} des racines de l'unité.\smallskip

Écrivons $r_K$ et $c_K$ les nombres de places réelles et complexes de $K$ et $\mathcal X_K=\Gal (M/K)$ le groupe de Galois de sa pro-$\ell$-extension abélienne $\ell$-ramifiée; $\E_K=\Zl\otimes _\ZZ E_K$ le $\ell$-adifié du groupe des unités $E_K$ de $K$ et $\,\U_{K_\ell}=\prod_{\l|\ell} \U_{K_\l}$ le produit des sous-groupes unités des $\R_{K_\l}$ aux places $\l$ qui divisent $\ell$. Soit enfin $s_K : \R_K\rightarrow \R_{K_\ell}=\prod_{\l|\ell} \R_{K_\l}$ l'épimorphisme naturel de semi-localisation.\smallskip

Avec ces notations, il vient (cf. e.g. \cite{Gi,Gra,Iw,J10,J17,J18,Mi2,Wa}):

\begin{Th}[{\bf Conjecture de Leopodt}]\label{CLG}
Pour tout corps de nombres $K$ et tout premier $\ell$, les conditions suivantes sont équivalentes:
\begin{itemize}
\item[(i)] Le corps $K$ admet (au plus) $c_K+1$ $\Zl$-extensions indépendantes.
\item[(ii)] Le groupe de Galois $\mathcal X_K=\Gal (M/K)$ est un $\Zl$-module de rang essentiel $c_K+1$.
\item[(iii)] L'application naturelle de $\;\E_K$ dans $\;\U_{K_\ell}$ est injective.
\item[(iv)] Les idèles principaux localement racines de l'unité sont les racines de l'unité: $\mu_K^{loc}=\mu^{\phantom{oc}}_K$.
\end{itemize}
Lorsque c'est le cas on dit que le corps $K$ satisfait la conjecture de Leopoldt  pour le premier $\ell$.
\end{Th}

\noindent{\em Preuve.} Ces équivalences sont bien connues (à l'exception de la dernière, qui l'est moins) et se démontrent immédiatement à l'aide de la Théorie $\ell$-adique du corps de classes: les $\Zl$-extensions étant $\ell$-ramifiées, les conditions {\em (i)} et {\em (ii)} sont évidemment équivalentes ; si $H$ désigne le $\ell$-corps de classes de Hilbert de $K$, i.e. sa $\ell$-extension abélienne non ramifiée maximale, il vient(cf. \cite{J31}):\smallskip

\centerline{ $\Gal(M/H) \simeq \U_K/s_K(\E_K)$,}\smallskip

\noindent où $s_K$ est le morphisme de semi-localisation, ce qui donne l'équivalence de {\em (ii)} et de {\em (iii)}, puisque le degré $[H:K]$ est fini et qu'on a: $\rg_{\Zl} (\U_K)=[K:\QQ]=r_K+2c_K$ et $\rg_{\Zl} (\E_K)=r_K+c_K-1$.

Enfin l'équivalence de {\em (iii)} et de {\em (iv)} résulte tout simplement du fait que le noyau de $s_K$ est formé d'idèles principaux qui sont localement partout des racines de l'unité.
\medskip

Ici encore le résultat d'indépendance de logarithmes de Baker-Brumer assure que cette conjecture est vérifiée par les corps absolument abéliens et quelques autres. Dans le cas général, il est commode d'en introduire introduire le défaut $\delta_K^{\mathscr L}$, qu'on peut définir comme suit:

\begin{Sco}
Le groupe $\mu_K^{loc}=\R_K\cap\prod_\p\mu^{\phantom{loc}}_{K_\p}$ des idèles principaux qui sont localement partout racines de l'unité est le produit:
\qquad\qquad$\mu_K^{loc}= \mu^{\phantom{loc}}_K  \Zl^{\delta^{\mathscr L}_K}$\smallskip

\noindent du sous-groupe des racines globales de l'unité et d'un $\Zl$-module libre de rang fini. 

On dit que $\delta_K^{\mathscr L}$ mesure le défaut de la conjecture de Leopoldt dans le corps $K$ pour le premier $\ell$.
\end{Sco}

\noindent{\em Preuve.} Comme sous-module de $\E_K\simeq \mu^{\phantom{loc}}_K \Zl^{r_K+c_K-1}$, le $\Zl$-module $\mu_K^{loc}$ est, en effet, le produit direct de son sous-module de torsion $\mu^{\phantom{loc}}_K$ et d'un $\Zl$-module libre de rang fini.

\Remarque Avec ces notations, le rang essentiel du $\Zl$-module $\mathcal X_K$ est ainsi: $c_K+1+\delta_K^{\mathscr L}$; et on retrouve ainsi les invariants classiques dans l'étude de ces questions.\medskip

Le groupe $\mu_K^{loc}$  est  un sous-groupe du groupe des unités logarithmiques $\,\wE_K$. L'image de celui-ci par l'application de semi-localisation $s_K:\;\R_K \rightarrow  \prod_{\l|\ell} \R_{K_\l}$ est ainsi un $\Zl$-module dont le rang essentiel fait apparaître la différence des défauts des conjectures de Gross-Kuz'min et de Leopoldt:

\begin{Cor}
On a: $\dim^{\phantom{c}}_{\Zl} s_K(\wE_K)=\dim^{\phantom{c}}_{\Zl}\wE_K -\dim^{\phantom{c}}_{\Zl}\mu_K^{loc}=r_K+c_K+\delta_K^{\mathscr G}-\delta_K^{\mathscr L}$.
\end{Cor}

\newpage
\noindent{\large \bf 5. Conditions algébriques suffisantes de la conjecture de Leopoldt}
\medskip

La trivialité du $\ell$-groupe des classes logarithmiques $\wCl_K$ implique trivialement la validité de la conjecture de Gross-Kuz'min pour le corps $K$ et le premier $\ell$, laquelle postule simplement la finitude du groupe  $\wCl_K$. En présence des racines $\ell$-ièmes de l'unité, elle emporte également, ce qui est moins évident, la validité de la conjecture de Leopoldt pour le corps $K$ et le premier $\ell$. Plus précisément:

\begin{Th}\label{CSL}
Soit $K$ un corps de nombres contenant les racines $2\ell$-ièmes de l'unité. Sous la condition $\wCl_K=1$, le corps $K$ satisfait les conjectures de Gross-Kuz'min et de Leopoldt pour le premier $\ell$ et il en est de même de chacun des étages $K_n$ de la tour cyclotomique $K_\infty/K$.
\end{Th}

\noindent{\em Preuve.} Que chacun des $K_n$ vérifie alors la conjecture de Gross-Kuz'min est immédiat: sous l'hypothèse $\wCl_K=1$, on a immédiatement $\T_K=1$, en vertu du lemme de Nakayama et, par suite,\smallskip

\centerline{ $\wCl_{K_n}=1$, pour tout $n\in\NN$,}\smallskip

\noindent puisque les $\ell$-groupes de classes logarithmiques des $K_n$ sont les quotients $\,\T_K/\T_K^{(\gamma^{\ell^n}-1)}$ relatifs aux sous-groupes $\Gamma^{\ell^n}$. Reste donc à vérifier que la trivialité de $\T_K$ entraîne la validité de la conjecture de Leopoldt à chaque étage fini $K_n$ de la tour cyclotomique.

L'hypothèse $\mu_{2\ell}^{\phantom{\ell}}\subset K$ assurant que $K_\infty$ contient toutes les racines $\ell$-primaires de l'unité, nous allons procéder par dualité en utilisant la Théorie de Kummer. Rappelons qu'en vertu des résultats de \cite{J18,J23}, l'application naturelle de $\;\mathfrak R_{K_n} = (\Ql/\Zl)\otimes_\ZZ K_n^\times$ dans $\;\mathfrak R_{K_\infty} = (\Ql/\Zl)\otimes_\ZZ K_\infty^\times$ est injective, ce qui permet de regarder les éléments de $\mathfrak R_{K_n}$ comme des radicaux kummériens au-dessus de $K_\infty$.

Considérons donc l'image canonique dans $\mathfrak R_{K_\infty}$ du sous-groupe $\mu^{loc}_{K_n}=\mathcal R_{K_n} \cap \prod_{\p_n}\mu_{K_{\p_n}}$ de $\mathcal R_{K_n}=\Zl\otimes_\ZZ K_n^\times$ formé des idèles principaux qui sont localement partout des racines de l'unité. Comme expliqué plus haut, la conjecture de Leopoldt pour le corps $K_n$ (et le premier $\ell$) affirme que le groupe 
$\mu^{loc}_{K_n}$ se réduit au sous-groupe $\mu_{K_n}$ des racines globales de l'unité. Lorsqu'elle est en défaut, le radical correspondant $(\Ql/\Zl)\otimes_{\Zl} \mu^{loc}_{K_n}$ est infini et définit une pro-$\ell$-extension abélienne non triviale de $K_\infty$ qui est, par construction, localement triviale partout, c'est-à-dire complètement décomposée en toutes les places; ce qui contredit la trivialité de $\T_K$. D'où le résultat.

\begin{Sco}
Soit  $K$ un corps de nombres contenant les racines $2\ell$-ièmes de l'unité. La conclusion du Théorème vaut encore sous la condition plus faible: $\T_K$ fini.
\end{Sco}

\noindent{\em Preuve.} D'un côté, la finitude du module de Tate $\T_K$ entraîne celle de tous ses quotients des genres $\wCl_{K_n}$, donc la conjecture de Gross-Kuz'min dans chacun des $K_n$. D'un autre côté, la preuve  montre que $\T_K$ est infini dès que la conjecture de Leopoldt est en défaut dans l'un des $K_n$.
\medskip

Plusieurs conditions suffisantes (mais non nécessaires) de la conjecture de Leopoldt pour un premier $\ell$ donné et un corps $K$ contenant les racines 2$\ell$-ièmes de l'unité, ont été proposées, notamment par Bertrandias et Payan \cite{BP}, en termes de plongements localement cyclotomiques; par Miki \cite{Mi1,Mi2}, en termes kummériens; par Gillard \cite{Gi}, en termes de théorie d'Iwasawa.

Comme expliqué dans \cite{J17,J18}, toutes ces conditions sont en fait rigoureusement équivalentes à la trivialité du groupe des classes logarithmiques, c'est-à-dire à la condition suffisante du Théorème \ref{CSL}. Par exemple la condition donnée dans \cite{Gi} s'écrit: $\T_K=1$. D'après le Théorème \ref{trivial}, elle s'écrit aussi bien $\;\wCl_K=1$.

L'intérêt majeur de la description que nous donnons, outre d'uniformiser ces diverses conditions, est de les traduire directement, via le corps de classes $\ell$-adique, en termes d'invariants arithmétiques effectivement calculables du corps considéré.\medskip

La notion de corps de nombres $\ell$-rationnel ou $\ell$-régulier (cf. \cite{GJ1,GJ2,JN,Mov,MN1}) et ses généralisations (cf. \cite{BJ,JS1,JS2}) permet ainsi de construire des tours infinies de corps de nombres non abéliens sur $\QQ$, qui satisfont les conjectures de Leopodt et de Gross-Kuz'min pour des raisons purement algébriques, et qui restent, pour l'instant, inaccessibles aux méthodes transcendantes.

\Remarque En présence des racines $2\ell$-ièmes de l'unité, transposée par la Théorie de Kummer, la conjecture de Leopodt exprime la finitude du sous-groupe des points fixes d'un tordu à la Tate convenable du module de Kuz'min-Tate $\,\T_K$. Ce point est précisé en fin de section 9.

\newpage
\noindent{\large \bf 6. Trivialité du sous-module ambige du module de Kuz'min-Tate $\mathcal T_K^{\,\Gamma}$}
\medskip

Plutôt que de tester la trivialité du module de Kuz'min-Tate $\mathcal T_K$ sur son quotient des genres ${}^\Gamma\mathcal T_K\simeq\wCl_K$, intéressons-nous à son sous-groupe des points fixes $\mathcal T_K^{\,\Gamma}$.\smallskip

D'évidence, on a les implications: ${}_\Gamma\T_K=1 \Leftrightarrow\T_K=1\Rightarrow\mathcal T_K^{\,\Gamma}=1$; mais il peut arriver que $\mathcal T_K^{\,\Gamma}$ soit trivial sans que $\mathcal T_K$ le soit.\smallskip

La condition $\mathcal T_K^{\,\Gamma}=1$ est ainsi plus faible que la condition suffisante étudiée dans la section précédente. Néanmoins, et contrairement à une conjecture de Seo (cf. \cite{Se1}, \S3), il peut arriver qu'elle soit en défaut. Pour voir cela, partons de l'interprétation classique du groupe $\mathcal T_K^{\,\Gamma}$:

\begin{Th}\label{Kuz'min}
Le sous-groupe ambige $\mathcal T_K^{\,\Gamma}$ du module de Kuz'min-Tate s'identifie au quotient $\,\wE_K/N_{k_\infty/K}(\wE_{K_\infty})$  du $\ell$-groupe des unités logarithmiques par le sous-groupe des normes cyclotomiques d'unités logarithmiques:\smallskip

\centerline{$\mathcal T_K^{\,\Gamma} \simeq \wE_K/N_{k_\infty/K}(\wE_{K_\infty})$.}
\end{Th}

\noindent{\em Preuve.} Ce résultat est essentiellement bien connu au vocabulaire près (cf. \cite{Kuz,Se1,Se2,Se3} ou encore \cite{Kat} dans un contexte plus général), et remonte au travail de Kuz'min cité ci-dessus (\cite{Kuz}, Prop. 7.5), mais certaines des démonstrations proposées ailleurs étant erronées (du fait d'une confusion malheureuse entre unités logarithmiques au sens naïf et au sens $\ell$-adique), nous en redonnons ci-dessous une preuve très simple basée sur l'arithmétique des classes ambiges.

Faisons choix d'un ensemble fini $S$ de places finies de $K$ contenant les places au-dessus de $\ell$ et suffisamment gros pour que le $\ell$ groupe $\Cl_K^S$ des $S$-classes d'idéaux de $K$ soit trivial; et appliquons les isomorphismes des classes ambiges de Chevalley (cf. \cite{Ch,J18}) aux $\ell$-groupes de $S$-classes des corps $K_n$. Puisque par hypothèse les $S$-classes d'idéaux ambiges sont triviales (du fait que $S$ contient à la fois les places ramifiées dans $K_n/K$ et un système de générateurs des classes étendues), il existe un isomorphisme canonique:\smallskip

\centerline{$\Cl_{K_n}^{S\;\Gamma_{\!n}}\simeq (\E_K^S \cap N_{K_n/K}(\R_{K_n})/N_{K_n/K}(\E_{K_n}^S)$}\smallskip

\noindent obtenu en envoyant la classe ambige d'un $S$-idéal $\mathfrak A$ sur l'image par la norme d'un générateur $\alpha$ du $S$-idéal principal $\mathfrak A^{(\gamma-1)}$. Par passage à la limite projective, il suit:\smallskip

\centerline{$\mathcal T^{S\,\Gamma}_K=\varprojlim \Cl_{K_n}^{S\;\Gamma}\simeq (\E_K^S \cap N_{K_\infty/K}(\R_{K_\infty})/N_{K_\infty/K}(\E_{K_\infty}^S)$.}\smallskip

{\em (i)} Par la théorie du corps de classes, la limite projective à gauche s'identifie au sous-groupe des points fixes du groupe de Galois $\mathcal T^S_K$ de la pro-$\ell$-extension abélienne maximale de $K_\infty$ qui est non ramifiée partout et complètement décomposée aux places au-dessus de $S$. Et comme la montée dans la tour cyclotomique $K_\infty/K$ a épuisé toute possibilité d'inertie aux places de $S$ qui ne divisent pas $\ell$, le groupe $\mathcal T^S_K$ coïncide de ce fait avec le module de Kuz'min-Tate $\mathcal T_K$.

{\em (ii)} Comme expliqué dans la section 1, le groupe de normes $\E_K^S \cap N_{K_\infty/K}(\R_{K_\infty})$ au numérateur à droite n'est autre que le groupe des unités logarithmiques du corps $K$.

{\em (iii)} Enfin, les places de $S$ étant {\em au sens logarithmique} presque totalement inertes dans la tour cyclotomique $K_\infty/K$, le groupe des normes de $S$-unités $N_{K_\infty/K}(\E_{K_\infty}^S)$ se réduit au sous-groupe des normes d'unités logarithmiques; d'où l'isomorphisme annoncé.

\begin{Cor}
La conjecture de Seo \og$\;\mathcal T_K^{\,\Gamma}=1$ pour tout corps de nombres $K$\fg { } est fausse.
\end{Cor}

\noindent{\em Contre-exemple.} Soit $d\equiv 2 \,[mod\,3]$ un entier naturel sans facteur carré et $k$ le corps quadratique réel $\QQ[\sqrt d]$. Le nombre premier $\ell=3$ étant inerte dans $K/\QQ$, le 3-groupe des unités logarithmiques de $k$ coïncide avec le 3-adifié $\E'_k=\ZZ_3\otimes_\ZZ E'_k$ du groupe des 3-unités. Or, pour $d=257$, l'unité fondamentale de $k$ n'est pas norme d'une 3-unité dans le premier étage $k[\cos (2\pi/9)]/k$ de la $\ZZ_3$-extension cyclotomique $k_\infty/k$. Il vient donc:\smallskip

\centerline{$|\,\mathcal T_k^{\,\Gamma}\,|=(\wE_k:N_{k_\infty/k}(\wE_{k_\infty}))\ge 3$;}\smallskip

\noindent ce qui contredit la conjecture de Seo (\cite{Se1}, \S3).\smallskip

\Remarque On peut imaginer que parmi les corps quadratiques réels $k=\QQ[\sqrt d]$ avec $d=2 \,[mod\,3]$, qui ont un 3-groupe des 3-classes non trivial, il en est une infinité qui contredisent la conjecture.

\newpage
\noindent{\large \bf 7. Lien avec les groupes de {\em normes universelles} au sens de  Kuz'min}
\medskip

Dans \cite{Kuz} Kuz'min définit les {\em normes universelles globales} associées à la $\Zl$-extension cyclotomique $K_\infty/K$ à partir des suites cohérentes $(\varepsilon_n)_{n\in\NN}\in\prod_{n\in\NN}\,\E'_{K_n}$ construites sur les $\ell$-adifiés $\,\E'_{K_n}=\Zl\otimes_\ZZ E'_{K_n}$ des groupes de $\ell$-unités attachés aux divers étages de la tour. C'est également le vocabulaire employé notamment par  \cite{Grt,Kat,MR,MN2} dans divers contextes.\smallskip

La condition de cohérence requise, $N_{K_m/K_n}(\varepsilon_m)=\varepsilon_n$ pour  $ m\ge n$, impose ainsi aux $\varepsilon_n$ non seulement d'être normes globales dans chaque étage fini $K_m/K_n$ (i.e. d'être des unités logarithmiques), mais d'être des normes d'unités logarithmiques dans chaque étage fini $K_m/K_n$ de la tour. De fait, il est facile de voir que les éléments initiaux de ces suites sont ceux de $N_{K_\infty/K}(\wE_{K_\infty)}$:

\begin{PDef}
Les termes initiaux des éléments de  la limite projective $\,\N=\varprojlim \E'_{K_n}$ des groupes de $\ell$-unités pour les applications normes sont les normes cyclotomiques des unités logarithmiques, i. e. les éléments $\,\E_K^{\,\nu}$ de l'intersection $N_{K_\infty/K}(\wE_{K_\infty})=\underset{n\in\NN}{\cap}N_{K_n/K}(\wE_{K_n})$.\smallskip

En d'autres termes, on a l'identité: $\varprojlim \E'_{K_n}=\varprojlim \wE_{K_n}=\underset{n\in\NN}{\prod}N_{K_\infty/K_n}(\wE_{K_\infty})=\underset{n\in\NN}{\prod}\,\E^{\,\nu}_{K_n}$.\smallskip

Nous disons que $\,\E^\nu_{K_n}$ est le (pro-)$\ell$-groupe des normes logarithmiques attaché à $K_n$.
\end{PDef}

 Donnons une preuve rapide de ce résultat de Kuz'min (cf. \cite{Grt} pour une approche fonctorielle). Les arguments utilisés dans la démonstration du Théorème \ref{Kuz'min} montrent que pour tout ensemble fini $S$ de places de $K$ contenant les places au-dessus de $\ell$, on a l'égalité:\smallskip

\centerline{$\underset{n\in\NN}{\cap}N_{K_n/K}(\E^S_{K_n})=\underset{n\in\NN}{\cap}N_{K_n/K}(\E'_{K_n})=\underset{n\in\NN}{\cap}N_{K_n/K}(\wE_{K_n}) = N_{K_\infty/K}(\wE_{K_\infty})=\E^{\,\nu}_K$.}\smallskip

\noindent Partons donc d'un élément $\varepsilon\in N_{K_\infty/K}(\wE_{K_\infty})$ et montrons que c'est le premier terme d'une suite cohérente $(\varepsilon_n)_{n\in\NN}$. 
Écrivons pour tout $n$ entier: $\varepsilon=N_{K_n/K}(\eta_n)$ avec $\eta_n\in\wE_{K_n}$; prenons pour $\varepsilon_1$ une valeur d'adhérence de la suite $(N_{K_n/K_1}(\eta_n))_{n\ge 1}$; 
et observons que nous avons par construction $\varepsilon=N_{K_1/K}(\varepsilon_1)$ et $\varepsilon_1\in N_{K_\infty/K_1}(\wE_{K_\infty})$; puis itérons le procédé. Nous pouvons construire ainsi de proche en proche les $\varepsilon_n$, ce qui donne bien la suite cohérente annoncée.\medskip

Contrairement aux groupes d'unités logarithmiques $\,\wE_{K_n}$, les groupes normiques $N_{K_\infty/K_n}(\wE_{K_\infty})$ {\em ne sont pas} caractérisés par leurs seules propriétés locales: leur définition fait intervenir de façon irréductible l'arithmétique globale de la $\Zl$-extension cyclotomique $K_\infty/K_n$. En revanche, ils possèdent les meilleures propriétés fonctorielles qu'on peut attendre (cf. \cite{Grt} ou \cite{Kuz}):

\begin{Th}\label{NormesU}
Soit $\Lambda = \Zl [[\gamma -1]]$  l'algèbre d'Iwasawa  attachée à un groupe procyclique $\Gamma = \gamma^{\Zl}=\Gal(K_\infty/K)$. Notons  $r_K$ et $c_K$ les nombres de places réelles et complexes du corps $K$. Alors:\smallskip
\begin{itemize}
\item[(i)] La limite projective $\,\N=\varprojlim \E'_{K_n}$ des groupes de $\ell$-unités est le produit direct du module procyclique (éventuellement trivial) $\varprojlim \mu^{\phantom{l}}_{K_n}$ et d'un $\Lambda$-module libre de dimension $r_K+c_K$.\smallskip
\item[(ii)] Pour chaque $n\in\NN$ le groupe normique $\,\E^\nu_{K_n}=N_{K_\infty/K_n}(\wE_{K_\infty})$ est le produit direct du $\ell$-groupe $\mu^{\phantom{l}}_{K_n}$ des racines de l'unité dans $K_n$ et d'un $\Zl$-module libre de dimension $(r_K+c_K)\ell^n$.
\end{itemize}
\end{Th}

Compte tenu de l'isomorphisme $\mathcal T_K^{\,\Gamma} \simeq \wE_K/N_{k_\infty/K}(\wE_{K_\infty})$ donné par le Théorème \ref{Kuz'min} et de l'interprétation de la conjecture de Gross-Kuz'min donnée par le Théorème \ref{CGK}, il suit:

\begin{Sco}
Le corps $K$ vérifie la conjecture de Gross-Kuz'min si et seulement si le sous groupe normique $\,\E_K^{\,\nu}=N_{k_\infty/K}(\wE_{K_\infty})$ est d'indice fini dans le groupe des unités logarithmiques $\,\wE_K$; auquel cas le groupe  $\,\wE_K$ est la racine dans $\R_K$ du $\Zl$-module $\,\E^{\,\nu}_K=N_{k_\infty/K}(\wE_{K_\infty})$.\par
De façon générale, le quotient $\,\wE_K/\sqrt{\E^{\,\nu\phantom{l}}_K}$ est un $\Zl$-module libre de dimension $\delta_K^{\mathscr G}$. Il représente donc le module de défaut de la conjecture de Gross-Kuz'min dans le corps $K$.
\end{Sco}

\begin{Sco}\label{FGamma}
Soit $\,\F_K$ le sous-module fini du $\Lambda$-module noethérien $\,\T_K$ et  $\,\F_K^{\,\Gamma}=\F_K\cap\T_K^{\,\Gamma}$ son sous-module des points fixes. Indépendamment de toute conjecture, on a l'isomorphisme canonique:\smallskip

\centerline{$\F_K^\Gamma \simeq \sqrt{N_{k_\infty/K}(\wE_{K_\infty})}/N_{k_\infty/K}(\wE_{K_\infty})=\sqrt{\E^{\,\nu\phantom{l}}_K}/\E^{\,\nu}_K$.}
\end{Sco}

\newpage
\noindent{\large \bf 8. Sous-module fini $\F_K$ du module de Kuz'min-Tate  et capitulation}
\medskip

La conjecture de Gross-Kuz'min postulant la finitude du sous-module des points fixes $\T_K^{\,\Gamma}$ du module de Kuz'min-Tate, il est naturel de s'intéresser au sous-module fini $\,\F_K$ de $\T_K$: comme $\Lambda$-module noethérien, le groupe $\T_K$ possède un plus grand sous-module fini, qui n'est autre que l'ensemble de ses éléments pseudo-nuls. 
En termes de classes logarithmiques, celui-ci a une double interprétation arithmétique très simple (cf. \cite{JM}, \S2.c et \cite{LMN}, Th. 1.4 pour les groupes de $\ell$-classes; \cite{Kuz}, Th. 1.2 pour les $\ell$-groupes  logarithmiques, au vocabulaire près):

\begin{Th}\label{Cap}
Les sous-groupes de capitulation logarithmiques $\wt{C\!ap}_{K_n}=\Ker\,(\wCl_{K_n}\rightarrow\wCl_{K_\infty}^{\,\Gamma_{\!n}})$ sont ultimement constants dans la tour cyclotomique $K_\infty/K$ et ont pour limite projective le sous-module fini $\;\F_K$ du module de Kuz'min-Tate $\T_K$. Sous la conjecture de Gross-Kuz'min dans $K_\infty$, un résultat analogue vaut pour les conoyaux de capitulation logarithmiques dans $\,\wCl_{K_\infty}=\varinjlim\,\wCl_{K_n}$:

\centerline{$\wt{C\!ap}_{K_n}=\Ker\,(\wCl_{K_n}\rightarrow\wCl_{K_\infty}^{\,\Gamma_{\!n}}) \underset{n \gg 0} {\simeq} \F_K \underset{n \gg 0} {\simeq} \wt{C\!ocap}_n=\Coker\,(\wCl_{K_n}\rightarrow\wCl{}_{K_\infty}^{\,\Gamma_{\!n}})$.}
\end{Th}

\noindent{\em Preuve.} Du point de vue de Théorie d'Iwasawa, la question de la capitulation se présente comme suit (cf. \cite{Kuz}, \S 1 et \cite{JM} \S2.c): on dispose d'un module noethérien et de torsion $T$ sur l'algèbre d'Iwasawa $\Lambda = \Zl [[\gamma -1]]$ attachée à un groupe procyclique $\Gamma = \gamma^{\Zl}$~; on note $\omega_n =\gamma^{\ell^n} -1$; et on s'intéresse aux noyaux des morphismes de transition $j_{n\mapsto m}:T_n \rightarrow T_m^{\Gamma_{\!n}}$ pour $m \gg n \gg 0$ induits par la multiplication par $\omega_m / \omega_n$ entre les quotients $T_n = T/\omega_n T$ et $T_m = T/\omega_m T$ pour $n$ et $m-n$ assez grands~; ce qui revient à considérer le noyau $\Ker_n = \Ker\,(T_n \rightarrow T_\infty)$ du morphisme d'extension à valeurs dans la limite inductive $T_\infty$ des $T_n$\smallskip

 Notant $F$ le plus grand sous-module fini de $T$, on obtient directement:\smallskip

\centerline{$\Ker_n = \{ t+\omega_n T \in T/ \omega_n T \, | \, (\omega_m / \omega_n)\ t  =
0 \} \underset{m \gg n \gg 0}{=} (F + \omega_n T ) / \omega_n T \underset{n \gg 0} {\simeq} F$,}\smallskip

\noindent puisque, pour $m \ge n \gg 0$, les $\frac{ \omega_m}{\omega_n}$ sont étrangers au polynôme caractéristique de $T$.\smallskip 

C'est le résultat classique sur la capitulation, qui ne nécessite aucune hypothèse supplémentaire sur le module $T$. Regardons maintenant les conoyaux. Des égalités \smallskip

\centerline{$T_m^{\Gamma_{\!n}} =   \{ t+\omega_m T \in T/ \omega_m T \, | \, \omega_n t \in  \omega_m T  \}$ et $j_{n\mapsto m}(T_n)=\frac{ \omega_m}{\omega_n} T/\omega_m T $,}\smallskip

\noindent  on tire d'autre part:\smallskip

\centerline{$\Coker_n = \{ t+\frac{ \omega_m}{\omega_n} T \in T/ \frac{ \omega_m}{\omega_n} T \, | \,  \omega_n t  \in \omega_m T \} 
=(T^{\Gamma_{\!n}}+\frac{ \omega_m}{\omega_n} T ) / \frac{ \omega_m}{\omega_n} T
\underset{ n\gg 0}{=} (F +\frac{ \omega_m}{\omega_n} T ) / \frac{ \omega_m}{\omega_n} T \underset{m \gg n}
\simeq F$,}\smallskip

\noindent dès lors que les $T^{\Gamma_{\!n}}$ sont finis.
Prenant $T=\T_K$ et observant qu'on a encore $\T_{K_n}=T$, nous obtenons le résultat annoncé d'après le Théorème \ref{CGK}, la condition de finitude requise exprimant précisément la validité de la conjecture de Gross-Kuz'min à chaque étage $K_n$ de la tour $K_\infty$.\medskip

On retrouve ainsi, mais de façon plus simple car les classes logarithmiques sont particulièrement bien adaptées à la co-descente galoisienne dans la $\Zl$-extension cyclotomique, l'analogue parfait des résultats classiques sur la capitulation pour les $\ell$-groupes de $\ell$-classes, avec en prime un résultat beaucoup moins connu sur les conoyaux des morphismes d'extension.

Lorsque, de plus, l'invariant {\it mu}\, d'Iwasawa du $\Lambda$-module $\T_K$ est nul, ce qui est le cas  si le corps $K$ est absolument abélien (d'après un résultat de Ferrero et Washington) et conjecturé pour tout corps de nombres (cf. \cite{Wa}), on peut également dire mieux: transposant dans le cadre logarithmique les résultats de \cite{GJ}, on voit de même que la capitulation $\,\wt{C\!ap}_{K_n} \simeq \F_K$ est alors un {\em facteur direct} du groupe $\,\wCl_{K_n}$ pour tout $n$ assez grand.\smallskip

\Remarque En présence d'une conjugaison complexe $\tau$, i.e. lorsque le corps $K$ est une extension quadratique totalement imaginaire d'un sous-corps $K^+$ totalement réel, les divers $\Zl$-modules étudiés plus haut sont munis canoniquement d'une action de $\Delta=\{1,\tau\}\simeq\Gal(K_\infty/K_\infty^+)\simeq\Gal(K/K^+)$. Pour $\ell$ impair ils s'écrivent comme sommes directes de leurs composantes {\em réelles} et {\em  imaginaires} à l'aide des idempotents $e_\pm=\frac{1}{2}(1\pm\tau)$ de l'algébre $\Zl[\Delta]$. Le Théorème \ref{NormesU} montre que le groupe normique $\,\E^{\,\nu}_K=N_{k_\infty/K}(\wE_{K_\infty})$ est réel, de même donc que son radical. On conclut alors du Scolie \ref{FGamma} que  la composante imaginaire du module fini $\,\F_K$ est triviale: $\F_K^\ominus=1$.\par
On retrouve par là le fait bien connu  que la composante imaginaire $\,\T_K^\ominus$ du module de Kuz'min-Tate ne possède pas de sous-module fini non nul (cf. e.g. \cite{FG,Gb,Iw}).

\newpage
\noindent{\large \bf 9. Lien avec les noyaux étales sauvages de la $K$-théorie}
\medskip

L'étude des noyaux de la $K$-théorie des corps de nombres via la cohomologie galoisienne remonte aux travaux de Tate \cite{Ta}. Nombre de contributions ont eu pour objet de préciser ces relations, en lien éventuellement avec les conjectures de Leopoldt ou de Gross-Kuz'min  (cf. e.g. \cite{AM1,AM2,J24,Kah,KM1,KM2,Ko1,Ko2,KM1,KC,Ng0,Ng1,Ng2,Sc}). L'appellation {\em noyau étale sauvage} est due à Nguyen Quang Do.\smallskip

Le point de vue logarithmique a été introduit dans \cite{J27}  pour les noyaux sauvages $\W_2(K)$ en liaison avec l'expression explicite  du symbole de Hilbert, puis développé dans \cite{JSo}, et finalement dans \cite{JM} pour les noyaux étales supérieurs $\W_{2i}(K)$. 
Rappelons qu'en vertu des résultats de Schneider \cite{Sc}, ceux-ci peuvent être décrits, en présence des racines d'ordre $\ell$-primaire de l'unité, comme les quotients des copoints fixes du $i$-ième tordu à la Tate du module de Kuz'min-Tate $\,\T_K$:\smallskip

\centerline{$\W_{2i}(K) \,\simeq\, {}^{\Gamma}(\Tl^{\otimes i}\otimes_{\Zl}\T_K)$,}\smallskip

\noindent où $\Tl^{\otimes i}$ désigne la $|i|$-ième puissance tensorielle du module de Tate $\Tl=\varprojlim\mu_{K_n}^{\phantom{loc}}$, pour $i \ge 0$; du module opposé $\lT = \Hom (\mu^{\phantom{loc}}_{\ell^\infty},\Ql/\Zl)$,  pour $i \le 0$; et $\W_0(K)=  {}^\Gamma \T_K\simeq\wCl_K$, par construction.\medskip

En l'absence éventuelle des racines de l'unité, raisonnant dans $K'=K[\zeta_{\ell}]$, on a d'après \cite{JM}:

\begin{Th} \label{WK}
Soient $\ell$ un nombre premier, $K$ un corps de nombres dont la tour cyclotomique $K[\zeta_{\ell^\infty}]/K$ est procyclique, $K'\supset K[\zeta_\ell]$ une extension abélienne de $K$ de groupe $\Delta=\Gal[K'/K]$ d'ordre $d$ étranger à $\ell$. Notons $|\mu^{\phantom{loc}}_{K'}|=\ell^m$. Pour tout $i \in \ZZ$ on a l'isomorphisme de $\Zl [\Delta ]$-modules:\smallskip

\centerline{${}^{\ell ^m}\!\W_{2i}(K')  \simeq \mu _{\ell ^m}^{\otimes i}\otimes _{\Zl} \wi \Cl _{K'}$}\smallskip

\noindent entre le quotient d'exposant $\ell ^m$ du noyau étale sauvage $\W_{2i}(K')$ et le tensorisé $i$ fois par $\mu_{\ell ^m}$ du $\ell$-groupe des classes logarithmiques du corps $K'$.\par
En particulier, si $\bar\omega$ désigne le caractère de l'action de $\Delta$ sur le module $\lT$, il suit:\smallskip

\centerline{${}^{\ell ^m}\!\W_{2i}(K) ={}^{\ell ^m}\!\W_{2i}(K')^{e_1} \simeq \mu _{\ell ^m}^{\otimes i}\otimes _{\Zl} \wi \Cl _{K'}^{e_{\bar\omega^i}}$}\smallskip

\noindent où $e_\varphi=\frac{1}{d}\sum_{\tau\in\Delta}\varphi(\tau^{-1})\tau$ désigne l'idempotent de l'algèbre $\Zl[\Delta]$ associé au caractère $\varphi$; de sorte que le $\ell$-rang comme la trivialité du noyau sauvage étale $\W_{2i}(K)$ se lisent sur la $\bar\omega^i$-composante du $\ell$-groupe des classes logarithmiques du corps $K'$.
\end{Th}

Donnons un exemple (le miroir de Scholz)  illustrant ce dernier résultat pour $\ell=3$ (cf. \cite{JM}):

\begin{Cor}\label{WK2}
Soit $k=\QQ[\sqrt d]$ un corps quadratique  distinct de $\QQ[\sqrt{-3}]$) et $k^*=\QQ[\sqrt{-3d}]$ son reflet dans l'involution du miroir. Alors pour chaque ${i>0}$, les quotients d'exposant 3 des noyaux sauvages étales de $k$ sont donnés à partir des 3-goupes de classes logarithmiques de $k$ et $k^*$ par:
\begin{center}$
{}^3\W_{2i}(k) \simeq \left\{ \begin{array}{ll}
\mu_3^{\otimes i} \otimes_{\ZZ_3} \wi \Cl_k \quad & pour \  i \ pair, \\
\mu_3^{\otimes i} \otimes_{\ZZ_3} \wi \Cl_{k^*} \quad & pour \ i \ impair.
\end{array} \right.
$\end{center}
\end{Cor}

\Remarque Pour $k = \QQ[\sqrt{-3}]$, il vient directement: ${}^3\W_{2i}(k) \simeq \mu_3^{\otimes i} \otimes_{\ZZ_3} \wi \Cl_k = 1$.
\medskip

Venons en maintenant au noyau du morphisme d'extension déjà étudié par \cite {AM1,AM2,Coa, Gb,J24,Kah,KM1,KM2,Ng1,Ng2, Va}. Les arguments conduisant au Théorème \ref{Cap} nous donnent semblablement:

\begin{Cor} Sous les hypothèses du Théorème \ref{WK}  et pour tout caractère $\ell$-adique irréductible $\varphi$ du groupe $\Delta$, la $\varphi$-composante $\,\F_{K'}^{e_\varphi}$ du sous-module fini de $\,\T_{K'}$ s'identifie pour $n\gg 0$ à la $\varphi$-partie\smallskip
\begin{itemize}
\item[(i)]  ${C\!ap'}_{K'_n}^{\,e_\varphi}$ du sous-groupe  de $\,\Cl'_{K'_n}$ formé des $\ell$-classes qui capitulent dans $\,\Cl'_{K'_\infty}$;

\item[(ii)]  $\wi{C\!ap}_{K'_n}^{\,e_\varphi}$ du sous-groupe  de $\,\wCl'_{K'_n}$ formé des classes logarithmiques qui capitulent dans $\wCl_{K'_\infty}$;

\item[(iii)] $\lT^{\otimes i}\otimes_{\Zl} CW\!K_{2i}^{e_{\varphi\omega^i}}$ du tensorisé  $\lT^{\otimes i}\otimes_{\Zl} \Ker (W\!K_2(K'_n) \rightarrow W\!K_2(K'_\infty))$, pour tout $i\ne 0$.

\end{itemize}
\end{Cor}

\Remarque  Ces résultats ne présupposent ni la conjecture de Gross-Kuz'min dans $K'$,  i.e.  la finitude de $\W_0(K') $; ni celle de Schneider, qui postule celle des  $\W_{2i}(K')$ pour $i<-1$; ni celle de Leopoldt, qui affirme dans ce contexte la finitude du quotient $\W_{-1}(K')={}^{\Gamma}(\lT\otimes_{\Zl}\T_K)$ ou, de façon équivalente, celle du sous-groupe $(\lT\otimes_{\Zl}\T_K)^\Gamma$ (cf. e.g. l'appendice de \cite{JM}).

\newpage
\noindent{\large \bf 10. Radicaux kummériens en lien avec les normes cyclotomiques}
\medskip

Lorsque le corps $K_\infty$ contient le groupe $\mu^{\phantom{lc}}_{\ell^\infty}\!$ des racines de l'unité d'ordre $\ell$-primaire, ses $\ell$-extensions sont caractérisées par leurs radicaux kummériens regardés dans le groupe:\smallskip

\centerline{$\mathfrak R_{K_\infty}=(\Ql/\Zl)\otimes_\ZZ K^\times_\infty\simeq(\Ql/\Zl)\otimes_{\Zl}\R_{K_\infty}$.}\smallskip

\noindent D'après \cite{J18,J23}, ces radicaux vérifient la théorie de Galois dans la tour $K_\infty/K$, ce qui s'écrit:\smallskip

\centerline{$\mathfrak R_{K_n}=(\Ql/\Zl)\otimes_\ZZ K^\times_n\simeq(\Ql/\Zl)\otimes_{\Zl}\R_{K_n}={\mathfrak R}_{K_\infty}^{\Gamma_{\!n}}$,}\smallskip

\noindent un résultat analogue valant pour les radicaux hilbertiens (\cite{J18}, Prop. I.2.2 \& I.2.18):\smallskip

\centerline{$\mathfrak H_{K_n} =\big\{\ell^{-k} \otimes x \in {\mathfrak R}_{K_n} \mid x \in \J_{K_n}^{\ell^k}\wU_{K_n} \big\}=\mathfrak H_{K_\infty}^{\Gamma_{\!n}}$\quad avec \quad $\mathfrak H_{K_\infty}=\cup_{n\in\NN}\;\mathfrak H_{K_n}$.}\smallskip

\noindent Le groupe  $\mathfrak H_{K_n}$ est le radical kummérien de {\em  l'extension de Bertrandias-Payan} du corps $K_n$, i.e. i.e. du compositum  $K_n^{bp}$ des pro-$\ell$-extensions abéliennes de $K_n$ qui sont {\em localement} $\Zl$-plongeables  (cf. \cite{J18,J23,Ng0}). Trois sous-groupes de $\mathfrak H_{K_\infty}$ nous intéressent plus particulièrement:\smallskip

\begin{itemize}
\item[(i)] Le noyau de Tate $\,\mathfrak N_{K_\infty}\!$ dans $\,\mathfrak R_{K_\infty}\!$ des symboles universels à valeurs dans $K_2(K_\infty)$:\smallskip

\qquad${\mathfrak N}_{K_\infty}= \big\{\ell^{-k}\otimes x \in {\mathfrak R}_{K_\infty} \mid \{\zeta_{\ell^k}, x\} = 1\; \text{ dans }\; K_2(K_\infty)  \big\}$;\smallskip

\item[(ii)] Le radical $\mathfrak Z_{K_\infty}$ de la réunion $Z_\infty=\cup_{n \in \NN}\, Z_n$ des composées des $\Zl$-extensions $Z_n$ des  $K_n$:\smallskip

\qquad${\mathfrak Z}_{K_\infty} =\Rad(Z_\infty/K_\infty)= \big\{\ell^{-k}\otimes x \in {\mathfrak R}_{K_\infty} \mid K_\infty[\sqrt[\ell^k]{x}] \subset Z_\infty \big\}$;\smallskip

\item[(iii)] Le radical $\we_{K_\infty}$ construit sur le groupe $\wE_{K_\infty}$ des unités logarithmiques dans $\R_{K_\infty}$:\smallskip

\qquad$\we_{K_\infty}= \big\{ \ell^{-k}\otimes x \in {\mathfrak R}_{K_\infty} \mid x \in \wE_{K_\infty}=\cup_{n \in \NN}\; \wE_{K_n} \big\} $.
\end{itemize}

\begin{Th}\label{CoKo}
Sous la condition  $\mu^{\phantom{lc}}_{\ell^\infty}\!\subset K_\infty$, on a l'égalité et les deux inclusions~:
\begin{itemize}
\item[(i) ] ${\mathfrak N}_{K_\infty} ={\mathfrak E}_{K_\infty}^{\,\nu}$ est l'image dans $\,{\mathfrak R}_{K_\infty}\!$ du groupe normique $\E_{K_\infty}^{\,\nu} = \cup_{n\in\NN} \;\E_{K_n}^{\,\nu}$.
\item[(ii)] $\,{\mathfrak N}_{K_\infty} \subset {\mathfrak Z}_{K_\infty}$, avec égalité sous la conjecture de Leopoldt dans $K_\infty$ (i.e. dans tous les $K_n$);
\item[(iii)] ${\mathfrak N}_{K_\infty} \subset \we_{K_\infty}$, avec égalité sous la conjecture de Gross-Kuz'min dans $K_\infty$ (i.e. idem).
\end{itemize}
\end{Th}

\Preuve L'égalité {\em (i)} résulte directement du fait que les duaux de Pontrjagin respectifs des deux modules sont tous deux des $\Lambda$-modules libres de même dimension $c_K$ (cf. \cite{Grt}, Th.  et \cite{MN2}, Prop. 1.7). D'un autre côté d'après \cite{J24}, Th. 8, les assertions  {\em (ii)} et  {\em (iii)} ne sont qu'une reformulation des résultats de Coates (\cite{Coa}, Th. 4) et de Kolster (\cite{Ko1}, Th. 2.6).\medskip

Ainsi, sous les deux conjectures de Gross-Kuz'min et de Leopoldt, les divers radicaux coïncident {\em en haut de la tour} et s'expriment directement à partir des groupes d'unités logarithmiques.
Cependant, comme expliqué dans \cite{J23,J24}, la descente impose une torsion de l'action galoisienne:

\begin{Th}\label{RadH}
Toujours sous la condition $\mu^{\phantom{lc}}_{\ell^\infty}\!\subset K_\infty$, on a les identités de descente:\smallskip
\begin{itemize}
\item[(i)] Le noyau $\Tl\otimes_{\Zl}\mathfrak{N}_{K_n}$ du morphisme de Tate $(\mu^{\phantom{l}}_\infty\!\otimes_\ZZ K^\times_n)^{\Gamma_{\!n}}\twoheadrightarrow K_2(K_n)$ est le sous-groupe divisible maximal du groupe des points fixes $ (\Tl\otimes\mathfrak H_{K_\infty}\! )^{\Gamma_{\!n}}\subset (\Tl\otimes\mathfrak R_{K_\infty}\! )^{\Gamma_{\!n}} \simeq (\mu^{\phantom{l}}_\infty\!\otimes_\ZZ K^\times)^{\Gamma_{\!n}}$.\smallskip

\item[(ii)] Le tensorisé $\,\we_{K_n}=(\Ql/\Zl)\otimes_{\Zl} \wE_{K_n}$ est le sous-module divisible maximal du radical hilbertien $\,\mathfrak H_{K_n} = \mathfrak H_{K_n}$; d'où l'isomorphisme: $\,\we_{K_\infty}^{\,\Gamma_{\!n}}/ \we_{K_n}\simeq\,\wt{C\!ap}_{K_n}\underset{n \gg 0} {\simeq} \,\F_K$.
\smallskip

\item[(iii)] Le radical $\,{\mathfrak Z}_{K_\infty} =\big\{\ell^{-k}\otimes x \in {\mathfrak R}_{K_\infty} \mid K_\infty[\sqrt[\ell^k]{x}] \subset Z_n \big\}$ attaché au compositum $Z_n$ des $\Zl$-extensions de $K_n$ est le sous-groupe divisible maximal du groupe $(\lT\otimes \mathfrak H_{K_\infty}\! )^{\Gamma_{\!n}}$ \smallskip
\end{itemize}

En résumé, on a donc les trois identités:\smallskip

\centerline{$\lT \otimes {\mathfrak Z}_{K_n} = ( \lT \otimes {\mathfrak Z}_{K_\infty} )^{\Gamma_{\!n}}_{\mathrm{div}};
\qquad \we_{K_n} = (\we_{K_\infty})^{\Gamma_{\!n}}_{\mathrm{div}};\qquad 
\Tl\otimes {\mathfrak N}_{K_n} = (\Tl\otimes{\mathfrak N}_{K_\infty}\! )^{\Gamma_{\!n}}_{\mathrm{div}}$;}\smallskip
\end{Th}

\Preuve L'assertion {\em (i)} provient de la description cohomologique par Tate du $K_2$ d'un corps de nombres (cf. \cite{Ta}); l'assertion {\em (iii)} de la Théorie de Kummer (cf. \cite{J23}). Quant à l'assertion {\em (ii)}, elle résulte de la suite exacte reliant noyau hilbertien et classes logarithmiques (cf \cite{J23}, 2\S c):\smallskip

\centerline{$1 \rightarrow \we_{K_n}  \rightarrow {\mathfrak H}_{K_n} \rightarrow \wCl{}_{K_n}^{\,\mathrm{tor}}  \rightarrow 1$,}\smallskip

\noindent jointe à la finitude du sous-groupe de torsion $\,\wCl{}_{K_n}^{\,\mathrm{tor}}$ de $\,\wCl{}_{K_n}$.  Passant à la limite inductive  et prenant les points fixes par $\Gamma_{\!n}$, on obtient ensuite  la suite exacte longue\smallskip

\centerline{$1 \rightarrow \we_{K_\infty}^{\,\Gamma_{\!n}}  \rightarrow {\mathfrak H}_{K_\infty}^{\,\Gamma_{\!n}}= {\mathfrak H}_{K_n}  \rightarrow (\wCl{}_{K_\infty}^{\,\mathrm{tor}})^{\Gamma_{\!n}}   \rightarrow \cdots$}\smallskip

\noindent d'où l'on déduit, par le lemme du serpent, que le quotient $\,\we_{K_\infty}^{\,\Gamma_{\!n}}/\we_{K_n}$ s'identifie au sous-groupe de capitulation $\,\wt{C\!ap}_{K_n}$; lequel est ultimement isomorphe à $\,\F_K$ d'après le Théorème \ref{Cap}.

\newpage
\noindent{\large \bf 11. Descente radicielle et radicaux initiaux}
\medskip

Du fait des différences introduites par la torsion à la Tate,  la détermination  des radicaux $\mathfrak Z_{K_n}$ et  $\mathfrak T_{K_n}$ n'est pas du tout immédiate (cf. e.g. \cite{Gr1,MN2}), en dehors du cas simple étudié par Bertrandias et Payan dans \cite{BP}, où le radical hilbertien convient:

\begin{Prop}\label{coïncidence}
 Supposons $\mu^{\phantom{lc}}_{2\ell}\subset K$ et notons $\ell^m$ l'ordre de $\mu^{\phantom{l}}_K$ avec $m\ge1$. Sous la condition $\wt{\Cl}_K=1$, les sous-groupes de $\ell^m$\!-torsion du noyau universel de Tate $\,\mathfrak N_K$, du radical $\,\mathfrak Z_K$ et du groupe divisible $\,\wt{\mathfrak E}_K$ coïncident avec le radical initial de l'extension de Bertrandias-Payan:\smallskip

\centerline{${}_{\ell^m}\!{\mathfrak H}_K ={}_{\ell^m}\!{\mathfrak N}_K = {}_{\ell^m}\!{\mathfrak Z}_K = {}_{\ell^m}\!\we_K 
= \big\{ \ell^{-m}\otimes x \in {\mathfrak R}_K \mid x \in \wE_K \big\}$.}
\end{Prop}

\noindent  {\it Preuve.} Sous la condition $\,\wt{\Cl}_K=1$, le corps $K$ vérifie, en effet, les conjectures de Gross-Kuz'min et de Leopodt (cf. Théorème \ref{CSL}). On a de plus l'égalité $\,\mathfrak H_K=\we_K$, en vertu de la suite exacte :\smallskip

\centerline{$1 \rightarrow \we_{K_n}  \rightarrow {\mathfrak H}_{K_n} \rightarrow \wCl{}_{K_n}^{\,\mathrm{tor}}= \wCl{}_{K_n}  \rightarrow 1$,}\smallskip

\noindent déjà rencontrée lors de preuve du Théorème \ref{RadH}. Il suit:\smallskip

\centerline{${}_{\ell^m}\!{\mathfrak H}_K = {}_{\ell^m}\!\we_K \simeq(\ell^{-m}\Zl/\Zl)^{c_K}$,}\smallskip

\noindent en vertu du Scolie \ref{Ulog}; d'où le résultat annoncé, puisque $\,\mathfrak H_K$ contient $\,\mathfrak Z_K$ comme $\,\mathfrak N_K$.\medskip

En fait, comme expliqué dans \cite{J24}, c'est le le module fini $\,\F_K$ (ou plutôt son dual $\,\F_K^*$)  qui représente l'obstruction à 
l'égalité des radicaux initiaux aux niveaux finis de la tour cyclotomique $K_\infty/K$.  Plus précisément, les résultats de \cite{Gb} et  \cite{J24} peuvent être complétés comme suit:

\begin{Th}\label{Final}
 Supposons toujours $\mu^{\phantom{lc}}_{\ell^\infty}\subset K_\infty$ et écrivons $\ell^{n+m}$ l'ordre de $\mu^{\phantom{l}}_{K_n}$ pour $n>1$. Alors, sous les conjectures de Leopoldt et de Gross-Kuz'min, il y a équivalence entre les assertions:\smallskip
\begin{itemize}
\item[(i)] $\,\F_K$ est trivial; autrement dit les morphismes $\,\wCl^{\phantom{l}}_{K_n} \rightarrow \wCl_{K_\infty}^{\,\Gamma_{\!n}}$ sont ultimement  bijectifs.

\item[(ii)] Les groupes $\,\T^{\Gamma_{\!n}}_K$ sont tous triviaux; autrement dit, on a les égalités $\,\wE_{K_n}=\E_{K_n}^{\,\nu}$.\smallskip

\item[(iii)] Pour $n\gg 0$, il y a égalité entre les sous-groupes de $\ell^n$\!-torsion des trois radicaux:

\centerline{${}^{\phantom{l}}_{\ell^n}\!{\mathfrak N}_{K_{n-m}} = {}^{\phantom{l}}_{\ell^n\!}{\mathfrak Z}_{K_{n-m}} = {}^{\phantom{l}}_{\ell^n}\!\we_{K_{n-m}} 
= \big\{ \ell^{-n}\otimes x \in {\mathfrak R}_{K_{n-m}} \mid x \in \wE_{K_{n-m}} \big\}$.}
\end{itemize}
\end{Th}

\noindent{\em Preuve.} 
 Le $\Lambda$-module $\,\F_K$ étant un $\ell$-groupe fini, on a trivialement: $\,\F_K=1\Leftrightarrow \F^{\Gamma_{\!n}}_K=1$. Sous la conjecture de Gross-Kuz'min, l'équivalence des deux premières assertions résulte donc de la double inclusion $\,\F^{\Gamma_{\!n}}_K \subset \T^{\Gamma_{\!n}}_K \subset \F_K$ et des interprétations données par les Théorèmes \ref{Kuz'min} et \ref{Cap}.\par

Notons $\bar n=n-m$.
Sous la condition $\,\F_K=1$, les radicaux $\we_{K_{\bar n}}$ construits sur les normes cyclotomiques vérifient la théorie de Galois dans la tour $K_\infty/K$ (cf. Th. \ref{RadH}). Cela étant, le groupe de Galois $\Gamma_{\!\bar n}$ agissant trivialement sur les racines $\ell^n$\!-ièmes de l'unité, on obtient, comme annoncé:\smallskip

\centerline{$v\!{\mathfrak N}_{K_{\bar n}} = {}^{\phantom{l}}_{\ell^n}\!{\mathfrak N}_{K_\infty}^{\,\Gamma_{\!{\bar n}}} 
= {}^{\phantom{l}}_{\ell^n}\! \we_{K_\infty}^{\,\Gamma_{\!{\bar n}}} = {}^{\phantom{l}}_{\ell^n}\!\we_{K_{\bar n}} 
= \big\{ \ell^{-n}\otimes x \in {\mathfrak R}_{K_{\bar n}} \mid x \in \wE_{K_{\bar n}}
 \big\}$.}\smallskip

\noindent Et le même résultat vaut pour les radicaux $ {}^{\phantom{l}}_{\ell^n\!}{\mathfrak Z}_{K_{\bar n}}$ sous la conjecture de Leopoldt.\smallskip

Réciproquement, on sait depuis Iwasawa \cite{Iw} que le dual de Kummer $\,\F_K^*=\Hom(\F_K,\mu_{\ell_\infty}^{\phantom{lc}})$ du groupe fini $\,\F_K$ s'identifie au conoyau $\Lambda^{c_K}/\mathcal Z_{K_\infty}$ du groupe $\,\mathcal Z_{K_\infty}=\Gal(Z_\infty/K_\infty)$ dans le $\Lambda$-module libre qu'il détermine. 
Sous les conjectures de Leopoldt et de Gross-Kuz'min les calculs qui précèdent  montrent que le radical ${}^{\phantom{l}}_{\ell^n}\!\we_{K_{\bar n}}$ est l'ensemble des éléments de $\ell^n$-torsion du sous-module divisible maximal du groupe des points fixes $\we_{K_\infty}^{\,\Gamma_{\!{\bar n}}}$. 
Dans la dualité de Kummer entre $\we_{K_\infty}$ et $\we_{K_\infty}^* \subset \Lambda^{c_K}$, l'orthogonal $({}^{\phantom{l}}_{\ell^n}\!\we_{K_{\bar n}})^\bot$ de ${}^{\phantom{l}}_{\ell^n}\!\we_{K_{\bar n}}$ est ainsi pour $n\gg 0$ (et en notations additives) le sous-module $\omega_{\bar n}\Lambda^{c_K} + \ell^n\we_{K_\infty}^*$ de $\we_{K_\infty}^*$.
Considérons donc le quotient $\Lambda^{c_K}/ ({}^{\phantom{l}}_{\ell^n}\!\we_{K_{\bar n}})^\bot$:
\begin{itemize}
\item d'un côté, puisque c'est un quotient de $\Lambda^{c_K}/ \omega_{\bar n}\Lambda^{c_K} \simeq \Zl^{c_{K_n}}$, sa décomposition comme produit de groupes cycliques fait intervenir au plus $c_{K_{\bar n}}=c_K\ell^{\bar n}$ facteurs;

\item d'un autre côté, son sous-module $\,\we_{K_{\bar n}}^*/({}^{\phantom{l}}_{\ell^n}\!\we^{\phantom{l}}_{K_{\bar n}})^\bot$, qui est le dual de Kummer de ${}^{\phantom{l}}_{\ell^n}\!\we_{K_{\bar n}}$ est, lui, un $\ZZ/\ell^n\ZZ$-module libre de dimension $c_{K_{\bar n}}$;
\end{itemize}
Il vient donc: $\Lambda^c/ ({}^{\phantom{l}}_{\ell^n}\!\we_{K_n})^\bot  \simeq  \sum_{i=1}^{c_{K_{\bar n}}} \, \ZZ/\ell^{n+e_i}\ZZ, \quad {\rm avec} \quad \prod_{i=1}^{c_{K_{\bar n}}} \ell^{e_i} = (\Lambda^{c_K}:\we_{K_{\bar n}}^*) =\;\mid \F_K \mid.$

Or l'égalité de ${}^{\phantom{l}}_{\ell^{n}}\!\mathfrak N_{K_{\bar n}}$ avec ${}^{\phantom{l}}_{\ell^{n}}\!\we_{K_{\bar n}}$ signifie que $({}^{\phantom{l}}_{\ell^n}\!\we_{K_{\bar n}})^\bot$ reste inchangé lorsqu'on tord l'action de $\Gamma_{\!\bar n}$  par le caractère cyclotomique, i.e. que $1+\ell^n$ agit trivialement sur le quotient $\Lambda^{c_K}/ ({}^{\phantom{l}}_{\ell^n}\!\we_{K_{\bar n}})^\bot$, qui est donc d'exposant $\ell^n$; ce qui donne: $\,\F_K=1$. De même pour  ${}_{\ell^{n}}\!\mathfrak Z_{K_{\bar n}}={}_{\ell^{n}}\!\we_{K_{\bar n}}$.

\Remarque
Sous les conjectures de Schneider \cite{Sc}, la condition $\,\F_K=1$ entraîne de façon semblable l'égalité des radicaux initiaux associés aux divers tordus à la Tate ${}_{\ell^n}^{\phantom{l}}\!\big(\lT^{\otimes i}\otimes(\Tl^{\otimes i}\otimes\mathfrak N_{K_\infty})^{\Gamma_{\bar n}}_{div}\big)$.

\newpage
\centerline{\sc Appendice}\medskip

\noindent{\large \bf 12. Le principe de Hasse dans la Théorie $\ell$-adique du corps de classes}
\medskip

Classiquement le principe de Hasse affirme que dans une extension cyclique $L/K$ de corps de nombres les éléments de $K^\times$ qui sont normes locales partout sont exactement les normes globales.\smallskip

Par plongement du groupe multiplicatif $K^\times$ dans le groupe des idèles $J_K$,  cela s'écrit:

\begin{Th*}[Principe de Hasse]
Dans une extension cyclique $L/K$ de corps de nombres, on a:\smallskip

\centerline{ $K^\times\cap N_{L/K}(J_L) = N_{L/K}(L^\times)$.}
\end{Th*}

Et ce résultat s'étend trivialement au cas procyclique comme suit:

\begin{Sco*}
Dans une extension procyclique $L$ d'un corps de nombres $K$, on a de même:\smallskip

\centerline{ $K^\times\cap N_{L/K}(J_L) = N_{L/K}(L^\times)$.}\smallskip

\noindent en convenant de définir les sous-groupes normiques global et local comme intersections respectives des mêmes sous-groupes associés aux sous-extensions (cycliques) finies $F/K$ de $L/K$:\smallskip

\centerline{$ N_{L/K}(L^\times)= \underset{K\subset F\subset L}{\cap} N_{F/K}(F^\times)\quad\&\quad  N_{L/K}(J_L)= \underset{K\subset F\subset L}{\cap} N_{F/K}(J_F)$.}
\end{Sco*}

Le but du présent appendice est d'expliquer comment ce principe se transpose dans le contexte de la Théorie $\ell$-adique du corps de classes. Bien qu'utilisé ici et là dans la description des unités logarithmiques, ce résultat ne fait pas l'objet, en effet, d'un énoncé explicite dans les textes fondateurs de l'arithmétique $\ell$-adique comme \cite{J18} ou \cite{J31}. Il nous a donc paru utile d'en donner ici une courte preuve pour la commodité du lecteur, les groupes de normes cyclotomiques jouant un rôle central dans le présent article.\smallskip

Rappelons que si $\R _{K_\p}=\varprojlim K_\p ^\times /K_\p^{\times \ell^m}$ désigne le compactifié $\ell$-adique du groupe multiplicatif $K_\p^\times$ du  complété de $K$ en la place $\p$, le {\em $\ell$-adifié  du groupe  des idèles} est défini comme le le produit $\J_K=\prod_\p ^{res}\R_{K_\p}$ des compactifiés $\R _{K_\p}$ restreint aux familles $(\x_\p)_\p$ dont presque tous les éléments tombent dans le sous-groupe unité  $\,\U_K=\prod_\p\U_{K_\p}$, où $\,\U_{K_\p}$ s'identifie au $\ell$-groupe $\mu_{K_\p}$ des racines de l'unité dans $K_\p$ d'ordre $\ell$-primaire, pour $\p\nmid\ell$; à celui des unités principales, pour $\p\mid\ell$. Et que le plongement diagonal de $K^\times$ dans le groupe des idèles $J_K$ induit un morphisme injectif du $\ell$-adifié $\,\R_K=\Zl\otimes K^\times$ du groupe multiplicatif $K^\times$ dans le $\ell$-adifié $\,\J_K$; ce qui permet de regarder $\,\R_K$ comme un sous-module de $\,\J_K$ (cf. \cite{J18} ou \cite{J31}). Avec ces notations, il vient:

\begin{Th}[Principe $\ell$-adique de Hasse]\label{PH}
Dans une $\ell$-extension cyclique de corps de nombres $L/K$, on a l'identité normique:

\centerline{ $\R_K \cap N_{L/K}(\J_L) = N_{L/K}(\R_L)$.}\smallskip

\noindent Et le même résultat vaut encore lorsque $L/K=K_\infty/K$ est une $\Zl$-extension, avec les conventions:\smallskip

\centerline{$ N_{K_\infty/K}(\R_{K_\infty})= \underset{n\in\NN}{\cap} N_{K_n/K}(\R_{K_n})\quad\&\quad  N_{K_\infty/K}(\J_{K_\infty})= \underset{n\in\NN}{\cap} N_{K_n/K}(\J_{K_n})$,}\smallskip

\noindent où $K_n/K$ décrit les étages finis de degrés respectifs $[K_n:K]=\ell^n$ dans la tour $K_\infty/K$.
\end{Th}

\noindent{\em Preuve.} Considérons une $\ell$-extension cyclique $L/K$ de degré, disons, $\ell^n$ et partons d'un idèle principal $\x^{\phantom{l}}_K\in\R_K\cap N_{L/K}(\J_L)$. Notons $x\mapsto x^\otimes$ le $\ZZ$-morphisme naturel de $K^\times$ dans $\,\R_K$ et observons que la classe de $\x_K$ modulo $\,\R_K^{\ell^n}$ est représentée par l'image $x_K^\otimes$ d'un élément
 $x^{\phantom{l}}_K$ de $K^\times$ (qui est défini de façon unique à une puissance $\ell^n$-ième près). Par construction, nous avons $\x^{\phantom{l}}_K=x_K^\otimes \y_K^{\ell^n}$ pour un $\y^{\phantom{l}}_K\in\R_K$. Par hypothèse, l'élément $x_K^\otimes$ est ainsi une norme locale i.e. contenu dans le sous-module $N_{L/K}(\J_L)$ de $\,\J_K$; de sorte que $x^{\phantom{l}}_K$ est lui-même dans le sous-groupe $N_{L/K}(J_L)$ de $J_K$. D'après le principe de Hasse classique, c'est donc la norme d'un élément global $y^{\phantom{l}}_L\in L^\times$. Il suit $\x^{\phantom{l}}_K=N_{L/K}(y_L^\otimes) \y_K^{\ell^n}=N_{L/K}(y_L^\otimes \y^{\phantom{l}}_K)\in N_{L/K}(\R_L)$, comme attendu.

\medskip

\noindent{\bf Remarque}. Pour chaque premier $p\neq\ell$, l'élément $p$ étant inversible dans l'anneau $\Zl$, chaque élément de $\,\R_K$ est une puissance $p^n$-ième pour tout $n\in\NN$, donc banalement norme dans toute $p$-extension $L/K$. Le principe $\ell$-adique de Hasse n'a donc d'intérêt que dans les (pro-)$\ell$-extensions.

\newpage
\centerline{\sc Index des principales notations}\medskip
\bigskip

Nous recensons ci-dessous les principales notations utilisées dans le corps de l'article.\bigskip

\noindent{\bf Notations latines}
\smallskip

$K$  : un corps de nombres arbitraire; $K_\p$: le complété de $K$ en la place $\p$;\
	
$K_{\infty} = \cup_{n \in \mathbb N} K_n$ avec $[K_n :K]= \ell^n$: la $\Zl$-extension cyclotomique de $K$;\

$K_n^{bp}$ : la composée des $\ell$-extensions abéliennes de $K_n$ localement $\mathbb Z_\ell$-plongeables;\

$Z_n$: le compositum des $\mathbb Z_\ell$-extensions de $K_n$; et $Z_\infty$ leur réunion;\

$K_\infty^{cd}$ : la pro-$\ell$-extension abélienne localement triviale maximale de $K_\infty$;\

$E_{K_n}$ : le groupe des unités de $K_n$ et $E'_{K_n}$ le groupe des $\ell$-unités;\

$r_K, c_K$: les nombres respectifs de places réelles et complexes du corps $K$;\

$\mathbb T_\ell = \varprojlim \mu_{\ell^n}$ : le module de Tate; et $\ov {\mathbb T}_\ell$: le module opposé (cf. \S 9).\

\medskip
\noindent{\bf Notations anglaises} (cf. \S 1)
\smallskip

$\R_K=\Zl\otimes_\ZZ K^\times$: le $\ell$-adifié du groupe multiplicatif du corps $K$;\

$\R_{K_\p} = \varprojlim K^\times_\p/K^{\times \ell^m}_\p\!\!$: le compactifié $\ell$-adique du groupe multiplicatif $K^\times_\p$;\

$\U_{K_\p}$: le sous-groupe unité et $\wU_{K_\p}$ le groupe des normes cyclotomiques dans $\R_{K_\p}$;\

$\J_K = \prod^{res}_\p \R_{K_\p}$ : le $\ell$-adifié du groupe des idèles de $K$;\

$\U_K=\prod_\p \U_{K_\p}$: le sous-groupe unité de $\J_K$;\

$\wU_K=\prod_\p \wU_{K_\p}$: le sous-groupe des normes cyclotomiques dans $\J_K$;\

$\E'_K=\Zl\otimes_\ZZ \E'_K$: le $\ell$-adifié du groupe des $\ell$-unités de $K$;\

$\E_K=\Zl\otimes_\ZZ \E_K=\R_K\cap\,\U_K$: le $\ell$-adifié du groupe des unités de $K$;\

$\wE_K=\R_K\cap\,\wU_K$: le $\ell$-groupe des unités logarithmiques de $K$;\

$\wE_{K_\infty}=\cup_{n\in\NN}\,\wE_{K_n}$: le $\ell$-groupe des unités logarithmiques de $K_\infty$;

$\E_K^{\,\nu}= N_{K_\infty/K}(\wE_{K_\infty})$: le sous-module des normes logarithmiques (cf. \S 7);

$\wCl_K=\J_K/\wU_K\R_K$: le $\ell$-groupe des classes logarithmiques du corps $K$;

$\T_K=\Gal(K^{cd}_\infty/K_\infty)$: le module de Kuz'min-Tate attaché au corps $K$ (cf. \S 2);

$\F_K$: le plus grand sous-$\Lambda$-module fini de $\,\T_K$ (cf. \S 8).

\medskip
\noindent{\bf Notations gothiques} (cf. \S 10)
\smallskip

${\mathfrak R}_{K_\infty} = (\mathbb Q_\ell/\mathbb Z_\ell) \otimes_\ZZ K_\infty^\times$: le radical universel attaché au corps $K_\infty$~;\par

${\mathfrak N}_{K_\infty}=\big\{\ell^{-m}\otimes x \in\mathfrak R_{K_\infty}\,|\, \{\zeta_{\ell^ m},x\}=1\big\}$: le noyau des symboles universels;

${\mathfrak H}_{K_\infty}=\big\{\ell^{-m}\otimes x \in\mathfrak R_{K_\infty}\,|\, (\frac{\zeta_{\ell^{{\scriptscriptstyle{\mathit m}}}},\,x}{\p_{\!\infty}})=1,\;\forall \p_\infty\big\}$: le noyau des symboles de Hilbert;

${\mathfrak Z}_{K_\infty}$: le radical du compositum des $\Zl$-extensions de tous les $K_n$;\par

$\mathfrak R_{K_n} = (\mathbb Q_\ell/\mathbb Z_\ell) \otimes_\ZZ K_n^\times=\mathfrak R_\infty^{\Gamma_{\!n}}$: le radical attaché au sous-corps $K_n$;

$\mathfrak N_{K_n}=\lT\otimes (\Tl\otimes{\mathfrak N}_{K_\infty}\! )^{\Gamma_{\!n}}_{\mathrm{div}}$: le noyau universel de Tate attaché au corps $K_n$.

${\mathfrak H}_{K_n}=\mathfrak H_\infty^{\Gamma_{\!n}}$: le radical hilbertien attaché au sous-corps $K_n$;\par

${\mathfrak Z}_{K_n}=\Tl\otimes (\lT\otimes{\mathfrak Z}_{K_\infty}\! )^{\Gamma_{\!n}}_{\mathrm{div}}$: le radical du compositum des $\Zl$-extensions de $K_n$;\par

${\mathfrak E}_{K_n}'$: le tensorisé $(\mathbb Q_\ell/\mathbb Z_\ell)\otimes_{\mathbb Z}E'_{K_n}$ du groupe des $\l$-unités;\par

$\we_{K_n}$: le tensorisé $(\mathbb Q_\ell/\mathbb Z_\ell\otimes_{\mathbb Z_\ell}\wE_{K_n}$ du groupe des unités logarithmiques;\par

$\mathfrak E^{\,\nu}_{K_n}$: l'image canonique de $\,\E_{K_n}^{\,\nu}$ dans $\,\mathfrak R_{K_n}$.\par

\medskip
\noindent{\bf Notations grecques}
\smallskip

$\mu^{\phantom{l}}_{\ell^n}$: le groupe des racines $\ell^n$-ièmes de l'unité et $\mu^{\phantom{l}}_{\ell^\infty\!}=\cup_{n\in\NN}\,\mu^{\phantom{l}}_{\ell^n}$;\

$\mu^{\phantom{l}}_K$, $\mu^{\phantom{l}}_{K_\p}$ : les $\ell$-groupes de racines de l'unité respectifs de $K$ et $K_\p$;\

$\mu^{loc}_K=\R_K\cap\,\prod_\p \mu^{\phantom{l}}_{K_\p}$, le groupe global des racines locales de l'unité;\

$\Gamma \ $: le groupe procyclique $\Gal(K_\infty/K)$ et $\gamma$ un générateur topologique de $\Gamma$ ;\

$\Lambda=\Zl[[\gamma-1]]$ l'algèbre d'Iwasawa du groupe $\Gamma$;\

$\Gamma_{\!n}$: le sous-groupe $ \Gamma^{\ell^n}=\Gal(K_\infty/K_n)$ et $\omega_n=\gamma^{\ell^n}-1$;\

$\delta^{\mathscr L}_K=\dim_{\Zl}\mu^{loc}_K/\mu^{\phantom{l}}_K$: le défaut de la conjecture deLeopoldt dans $K$;\

$\delta^{\mathscr G}_K=\dim_{\Zl}\wE_K/\E^{\,\nu}_K$: le défaut de la conjecture de Gross-Kuz'min dans $K$.\

\newpage

\noindent{\sc Remerciements}\smallskip

L'auteur remercie particulièrement Bill {\sc Allombert} qui a effectué à l'aide du logiciel {\sc pari} les calculs sur les classes et unités logarihmiques présentés dans ce travail, ainsi que le rapporteur pour les nombreuses et judicieuses observations contenues dans son rapport d'arbitrage.


\def\refname{\normalsize{\sc  Références}}

\bigskip\noindent
{\small
\begin{tabular}{l}
{Jean-François {\sc Jaulent}}\\
Institut de Mathématiques de Bordeaux \\
Université de {\sc Bordeaux} \& CNRS \\
351, cours de la libération\\
F-33405 {\sc Talence} Cedex\\
courriel : Jean-Francois.Jaulent@math.u-bordeaux1.fr 
\end{tabular}
}

 \end{document}